\documentclass[11pt]{amsart}

\usepackage[utf8]{inputenc}
\usepackage[T1]{fontenc}
\usepackage[english]{babel}

\usepackage{amssymb}
\usepackage{mathrsfs}
\usepackage{amsfonts}
\usepackage{amssymb}
\usepackage{soul, ulem}

\usepackage{amsmath,amsthm,bbm}

\usepackage{graphicx}

\usepackage{hyperref}

\usepackage{mathtools}
\usepackage[x11names]{xcolor}
\newtagform{blue}{\color{blue}(}{)}

\numberwithin{equation}{section}
\usepackage[nodayofweek]{datetime}

\newcommand{\dz}{\mathrm{d}z}

\def\dsum{\mathop{\displaystyle \sum }}%
\def\dprod{\mathop{\displaystyle \prod }}

\newcommand{\1}{\mathbbm{1}}

\def\pf{{\bf Proof: }}

\def\om{\omega}
\def\Om{\Omega}

\def\un{\infty}
\def\eh{\frac{1}{2}}

\def\epf{$\Box$ \\}

\def\bx{\overline{x}}
\def\by{\overline{y}}

\def\bxi{\overline{\xi}}
\def\beeta{\overline{\eta}}

\newtheorem{thm}{Theorem}[section]
\newtheorem{pr}[thm]{Proposition}
\newtheorem{co}[thm]{Corollary}
\newtheorem{lem}[thm]{Lemma}
\theoremstyle{definition}

\newtheorem{rem}[thm]{Remark}

\newcommand{\N}{\mathbb{N}}
\newcommand{\R}{\mathbb{R}}

\newcommand{\Z}{\mathbb{Z}}

\newcommand{\be}{\begin{equation}}
	\newcommand{\ee}{\end{equation}}
\newcommand{\bdm}{\begin{displaymath}}
	\newcommand{\edm}{\end{displaymath}}
\newcommand{\bean}{\begin{eqnarray}}
	\newcommand{\eean}{\end{eqnarray}}
\newcommand{\bea}{\begin{eqnarray*}}
	\newcommand{\eea}{\end{eqnarray*}}

\newcommand{\cB}{\mathcal{B}}

\newcommand{\cT}{\mathcal{T}}

\author[P. Imkeller]{Peter Imkeller}
\address[P. Imkeller]{Humboldt-Universit\"at zu Berlin, Institut f\"ur Mathematik, Unter den Linden 6, D-10099 Berlin, Germany}
\email{\tt imkeller@mathematik.hu-berlin.de}
\thanks{P.~Imkeller was supported in part by DFG Research Unit FOR 2402.}

\author[O. Menoukeu Pamen]{Olivier Menoukeu Pamen}
\address[O. Menoukeu Pamen]{Institute for Financial and Actuarial Mathematics, Department of Mathematical Sciences, University of Liverpool, L69 7ZL, United Kingdom}
\email{\tt menoukeu@liv.ac.uk}
\thanks{O.~Menoukeu Pamen acknowledges the funding provided by the Alexander von Humboldt Foundation, under the programme financed by the German Federal Ministry of Education and Research entitled German Research Chair No 01DG15010.}

\author[F. Proske]{Frank Proske}
\address[F. Proske]{Department of Mathematics, University of Oslo Ghana, Moltke Moes Vei 35, Oslo, P.O. Box 1053, Blindern, 0316, Norway}

\email{\tt proske@math.uio.no}

\title[Takagi type functions and existence of local time]{Takagi type functions and dynamical systems: the smoothness of the SBR measure and the existence and smoothness of local time}



\begin{document}
	\selectlanguage{english}

	\begin{abstract}
		
		We investigate the occupation measures and local times of Takagi-type functions with roughness parameter $\gamma$, which are H\"older continuous with exponent
		$H=\frac{\log\gamma}{\log(1/2)}.$ Analytical insight is obtained by embedding these functions into a dynamical system related to the baker transform, whose global attractor is the graph of the Takagi function. The associated stable manifolds support Sinai-Bowen-Ruelle (SBR) measures, which we identify with the laws of certain symmetric Bernoulli convolutions. Dually, where duality is induced by time reversal, we derive a representation of the Takagi-type curves centred around the stable fibres in terms of Bernoulli convolutions, thereby relating SBR measures to occupation measures. While Bernoulli convolutions belong to the first Rademacher chaos, we show that the occupation measure is naturally represented as a functional in the second Rademacher chaos within the framework of non-Gaussian Malliavin calculus.
		
		Using a Fourier-analytic criterion together with variants of Weyl's equidistribution theorem, we prove that Takagi-type curves admit square-integrable local times for $\gamma=2^{-1/m}, \, m\geq 7,$
		and that the same conclusion holds for drifted Takagi curves for almost every $\gamma\in(1/2,1)$. 

	\end{abstract}
	
	\maketitle

	\section{Introduction}

\begin{figure}[h]
		\includegraphics[height=12cm,width=16cm]{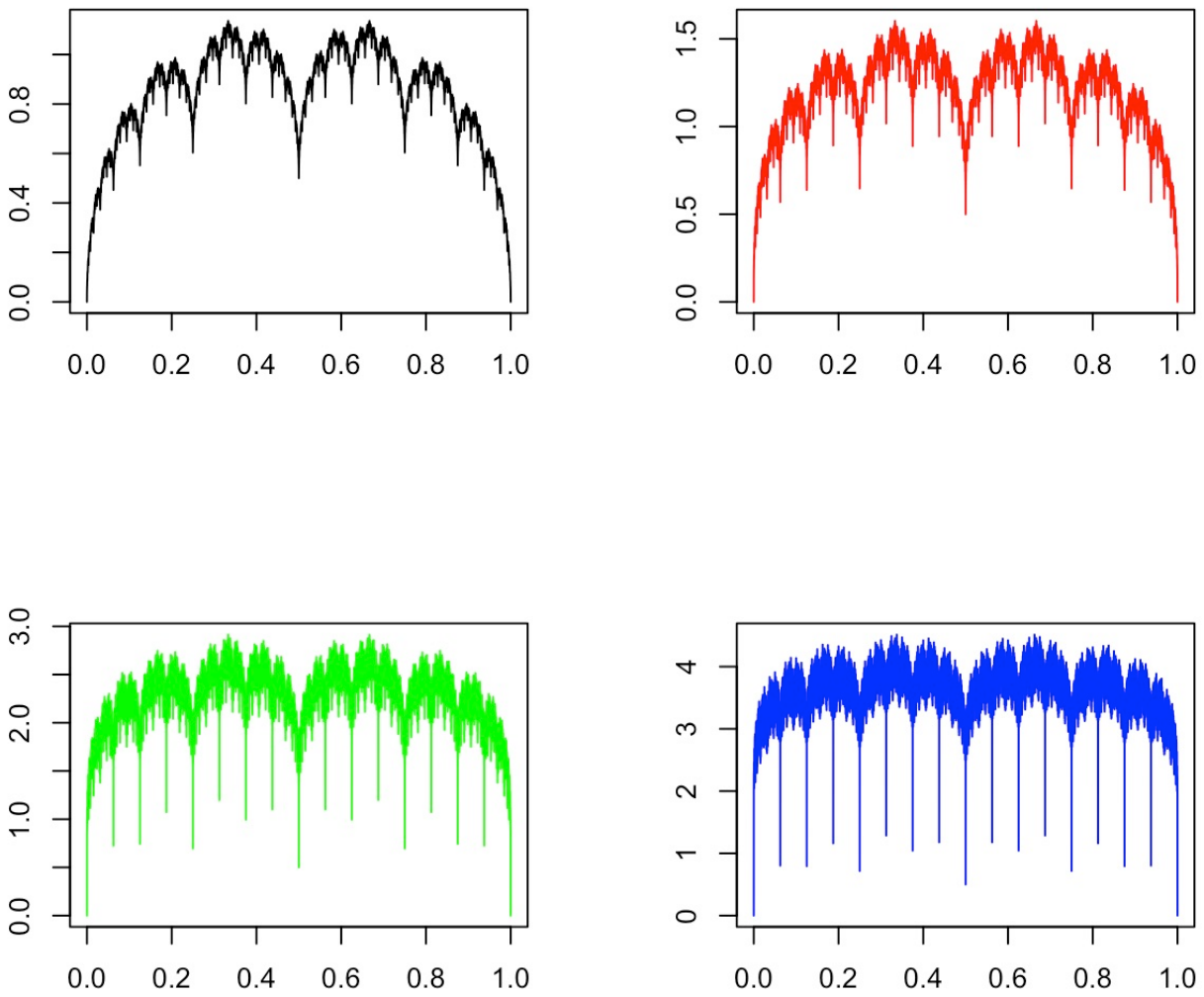}
		\caption{Takagi curve for $\gamma=2^{-\frac{1}{2}}$, \textcolor{red}{$\gamma=2^{-\frac{1}{3}}$}, \textcolor{green}{$\gamma=2^{-\frac{1}{6}}$}, \textcolor{blue}{$\gamma=2^{-\frac{1}{10}}$}}
		\label{Graph1}
	\end{figure}

	The interest in the rough Takagi-type curves studied in this paper arose from
	the study of a two-dimensional example of such functions in the context of the
	Fourier analytic approach to rough path analysis and rough integration
	developed in
	\cite{gubinelliimkellerperkowski2015,gubinelliimkellerperkowski2016}. In \cite{gubinelliimkellerperkowski2016}, the construction of a
	Stratonovich-type integral of a rough function $f$ with respect to another
	rough function $g$ is based on the notion of \emph{paracontrol}, a Fourier
	analytic extension of Gubinelli's notion of controlled rough paths
	\cite{gubinelli2004}. In search of natural examples of two-dimensional rough
	functions for which neither component is controlled by the other,
	Imkeller and Pr\"omel \cite{imkellerproemel15} introduced a pair of
	Weierstrass functions $W=(W_1,W_2)$. The two components oscillate on every dyadic scale with a phase shift of $\pi/2$: whenever one component has minimal increments, the other has maximal
	ones, and vice versa. This complementary oscillatory behaviour explains why
	the two components fail to control one another in the sense of Gubinelli. Moreover, the L\'evy areas of the finite approximations fail to converge, revealing a remarkable geometric pathology. 
	This naturally led us to investigate a closely related family of deterministic rough curves that
	exhibits comparable roughness but having a richer geometric structure
	amenable to the analysis developed in the present paper, namely the
	Takagi-type functions
	
	\[
	\mathcal T(x)=\sum_{n=0}^{\infty}\gamma^n\Phi(2^n x),
	\qquad x\in[0,1],
	\]
	where $\Phi(x)=d(x,\mathbb Z)$ denotes the distance from $x$ to the nearest
	integer and $\gamma\in(1/2,1)$ is a roughness parameter (see
	Figure~\ref{Graph1}). These curves are H\"older continuous with Hurst exponent
	\[
	H=\frac{\log\gamma}{\log(1/2)},
	\]
	and retain many of the geometric features of the Weierstrass functions while
	being considerably more amenable to the Fourier analytic methods developed
	below.
	
	We continue the study of the geometric properties of Takagi-type curves by
	investigating their occupation measures and local times. More precisely, we
	ask under which conditions on $\gamma$ the Takagi curve $\mathcal T$ admits a
	square-integrable local time. We first answer this question for a drifted
	Takagi curve $H$, introduced in \eqref{eq:FunctionH}, obtained by perturbing
	$\mathcal T$ by a smooth path arising naturally from the associated
	dynamical system, and subsequently for the classical Takagi function itself. These questions are closely related to the theory of regularization by noise.
	A fundamental example is Brownian motion, whose local time plays a central
	role in the path-by-path analysis of differential equations. Davie
	\cite{davie2007} proved path-by-path uniqueness for stochastic differential
	equations with bounded measurable drift by exploiting Brownian trajectories.
	Subsequently, Catellier and Gubinelli \cite{catelliergubinelli16} showed that
	paths possessing sufficiently regular local times have a remarkable
	regularizing effect on otherwise ill-posed ordinary differential equations.
	Their theory applies to a broad class of stochastic processes, including
	fractional Brownian motion. This initiated an active line of research on
	regularization by noise; see, for example,
	\cite{aminebanosproske2017,galeatigubinelli21,harangperkowski2021,priola2018,shaposhnikov2016}
	for ordinary differential equations,
	\cite{beflagumau19,butkovskymytnik2019} for partial differential equations,
	and
	\cite{bogsodieyepamen2022,bogsodieyepamen2023,bogsopamen2023}
	for two-parameter equations driven by the Brownian sheet. A central concept introduced in \cite{catelliergubinelli16} is
	$(\rho,\gamma)$-irregularity, which is characterised through Sobolev
	regularity of the occupation measure, or equivalently by the decay of its
	Fourier transform. It is shown there that perturbing an ill-posed equation by
	a $(\rho,\gamma)$-irregular path typically restores well-posedness. However,
	in all of the above works the perturbation is stochastic. As pointed out by
	Flandoli \cite{flandoli2015}, it is therefore natural to ask whether
	stochastic perturbations can be replaced by deterministic functions that are
	sufficiently ``strongly spread''. Since this spreading property is naturally
	quantified by occupation measures and their local times, Takagi-type curves
	provide a natural deterministic test case.
	
	Very little is known about occupation measures of deterministic rough curves. To the best of our knowledge, the only previous result in this direction is due to Buczolich \cite{Buczolich08}, who proved that, for the classical Takagi function corresponding to $\gamma=\frac12$, the occupation measure is singular with respect to Lebesgue measure. This naturally raises the question of what happens in the rough regime $\gamma>1/2$. The present paper answers this question by establishing the existence of square-integrable local times for a broad class of Takagi-type curves. The main conceptual contribution of the paper is the establishment of a transfer principle that allows results from the theory of Bernoulli convolutions to be carried over to occupation measures of deterministic Takagi-type curves. Exploiting a duality induced by the underlying hyperbolic dynamics, we show that the Fourier transform of the occupation measure is governed by the same probabilistic structures as Bernoulli convolutions. This enables fundamental results from the latter theory, including Solomyak's celebrated solution of the Erd\H{o}s problem, to be carried over to the setting of deterministic rough curves.
	
	Following the work of Hunt \cite{hunt1998}, Bara\'nski
	\cite{baranski02,baranski12,baranski14published,baranski15survey},
	Keller \cite{keller17-Publication-of-2015}, Shen
	\cite{shen2017-Publication-of-2015}, and others, we interpret the
	Takagi-type curves as global attractors of hyperbolic dynamical systems
	generated by suitable baker transformations. This dynamical systems
	viewpoint has proved extremely fruitful in the study of the geometric
	properties of Weierstrass-type functions, particularly for determining
	the Hausdorff dimension of their graphs (see the survey
	\cite{baranski14published}).
	
	For Takagi-type curves, we adopt the same framework. Besides the spatial
	variable $x$, which encodes the expanding direction of the baker
	transformation, we introduce an auxiliary variable $\xi$ describing the
	contracting direction. The forward dynamics simultaneously expands $x$
	by the factor $2$ and contracts $\xi$ by the factor $\frac12$, while
	these roles are interchanged under time reversal. This naturally leads
	to a three-dimensional hyperbolic dynamical system $F$ whose global
	attractor is the Takagi curve $\cT$. Associated with the smallest Lyapunov exponent of the linearisation of $F$ is a stable
	manifold generated by the Weierstrass-type series
	\[
	S(\xi,x)=
	-\sum_{n=1}^{\infty}
	\kappa^n\Phi'\!\left(B_2^n(\xi,x)\right),
	\]
	where $\kappa=\frac{1}{2\gamma}\in(1/2,1)$ is the roughness parameter
	dual to $\gamma$. The pushforward of Lebesgue measure by
	$S(\cdot,x)$ for $x\in [0,1]$-fixed coincides with the $x$-marginal of the
	Sinai-Bowen-Ruelle (SBR) measure of $F$. The definition of $F$ as a linear transformation combined with a smooth
	perturbation naturally reflects the self-affine structure of the Takagi
	curve. Self-affinity provides a mechanism for relating the microscopic
	geometry of the curve to its macroscopic behaviour through successive
	rescalings. Our main analytical tool for exploiting this structure is a
	system of \emph{telescoping relations}, inspired by the approach of Keller
	\cite{keller17-Publication-of-2015}. By developing these relations in both
	the forward and backward time directions, we uncover a duality, induced by
	time reversal, between the Sinai-Bowen-Ruelle measure and the occupation
	measure associated with the local time.
	
	More precisely, our analysis is centred around the doubly infinite series
	$$
	H(\xi,x)
	= \sum_{n\in\Z} \gamma^{-n} \left[\Phi \big( B^{-n}_2(\xi, x)\big) - \Phi \big( B^{-n}_2(\xi, 0)\big)\right],
	\quad \xi, x\in[0,1].
	$$
	A fundamental identity links $H$, the Takagi curve $\mathcal T$, and the
	stable manifold function $S$:
	$$
	H(\xi,y) - H(\xi,x) = \cT(y) - \cT(x) - \int_x^y S(\xi,z) \dz.
	$$
	This identity admits a natural geometric interpretation. If
	$l_{(\xi,x,w)}$ denotes the stable fibre through the point
	$(x,\mathcal T(x))$, defined by
	$$
	\frac{\mathrm{d}}{\mathrm{d} v}l_{(\xi,x,w)}(v)=S(\xi,v),
	\qquad
	l_{(\xi,x,w)}(x)=w,
	$$
	with $w=\mathcal T(x)$, then
	$$l_{(\xi, y, \cT(y))}(y) - l_{(\xi, x, \cT(x))}(y) = H(\xi,y) - H(\xi,x),\quad \xi, x, y\in[0,1].$$
	Thus, the increments of $H$ measure the vertical distance between stable
	fibres of the attractor.
	
	This identity reveals a remarkable duality. The Sinai-Bowen-Ruelle measure is naturally described by analysing the stable manifold function $S(\xi,\cdot)$  along the jumps of the contracting variable $\xi$. Dually, the occupation measure of the Takagi curve is obtained by analysing the increments of $H(\cdot,x)$  along the jumps of the expanding variable $x$. Thus, the roles of the expanding and contracting directions are interchanged through time reversal. To analyse the Sinai-Bowen-Ruelle measure, we decompose the stable manifold function $S(\xi,\cdot)$ along the dyadic jumps of $\xi$. This decomposition gives rise to a symmetric Bernoulli convolution, thereby revealing an unexpected connection between the geometry of Takagi-type curves and one of the classical objects of fractal analysis. The rich theory of Bernoulli convolutions, culminating in the metric regularity results of Erd\H{o}s, Solomyak and Wintner, then becomes available as an analytical tool. The main objective of the present paper is to show that this connection has genuine geometric consequences for deterministic Takagi-type rough curves by establishing regularity properties of their occupation measures.

	For the investigation of the occupation measure, in a dual step, we decompose
	$H(\cdot,x)$ along the jumps of the dyadic components of $x$ from $0$ to $1$.
	This leads to a surprising representation in which the essential part of
	$H(\cdot,x)$ is expressed as a functional involving a mixture of integrals in
	the first and second Rademacher chaoses. The Rademacher version of Malliavin
	calculus was developed in \cite{nourdinpeccatireinert10} and subsequent works.
	In this framework, the role of multiple Wiener integrals is played by multiple
	integrals with respect to a sequence $(X_n)_{n\in\mathbb N}$ of i.i.d.
	Bernoulli random variables. The classical Bernoulli convolution (and, in
	particular, the SBR function $S(\cdot,0)$) can be identified with a first
	Rademacher chaos integral
	\[
	I_1(f)=\sum_{n\in\mathbb N}f(n)X_n,
	\]
	where $f\in\ell_2(\mathbb N)$, whereas the second chaos is generated by
	double integrals of the form
	\[
	I_2(g)=\sum_{n,m\in\mathbb N}g(n,m)X_nX_m,
	\]
	where $g\in\ell_2(\mathbb N^2)$ vanishes on the diagonal. By computing the
	mixed moments
	\[
	E\!\left[I_1(f)^nI_2(g)^m\right], \qquad n,m\in\mathbb N,
	\]
	we derive in Proposition~\ref{II:FT-externalparts} a semi-explicit representation of the Fourier transform of the essential part of the occupation measure. This representation constitutes the principal analytical ingredient of the paper. It establishes a bridge between the Fourier transform of the occupation measure of the drifted Takagi curve and the characteristic function of an infinite Bernoulli convolution. Through this connection, Plancherel's theorem together with the $L^2$-regularity theory for Bernoulli convolutions yields our main result: for almost every $\gamma\in(1/2,1)$, the drifted Takagi curve admits a square-integrable local time.	More broadly, our results reveal that a regularity phenomenon from the theory of Bernoulli convolutions also governs the occupation measures of deterministic Takagi-type curves. While Bernoulli convolutions arise from a probabilistic construction, we show that analogous $L^2$-regularity emerges in a deterministic self-affine geometric setting. In this sense, our theorem may be viewed as a deterministic fractal counterpart of the metric regularity phenomenon established for Bernoulli convolutions.

	To establish the square integrability of the Fourier transform in the
	non-drifted case, we derive a product representation whose logarithm is an
	average of the periodic function $\log\cos^2$. Thus, the problem is
	reduced to understanding the equidistribution modulo $\pi$ of the arguments
	appearing in this product. This naturally leads to metric number theory,
	where we combine Koksma's refinement of Weyl's equidistribution theorem
	\cite{koksma35,weyl10} with quantitative \underline{discrepancy} estimates due to
	Erd\H{o}s and Tur\'an (see Kuipers and Niederreiter
	\cite{kuipersniederreiter74}).  Applying these techniques, we prove that the
	Fourier transform of the occupation measure is square integrable for
	$\gamma=2^{-1/m}$ with $m\ge7$. This parameter regime lies within the range
	$m\ge1$ for which Wintner \cite{wintner35} established the absolute
	continuity of the corresponding Bernoulli convolutions. As a consequence, we
	show that the classical Takagi function admits a square-integrable occupation
	density. To the best of our knowledge, this is the first result establishing
	the existence of square-integrable local times for the classical Takagi
	function and, more generally, deterministic Takagi-type rough curves.

	The paper is organized as follows. In Section~\ref{s:attractor}, we interpret Takagi-type curves as attractors of a dynamical system generated by the baker transform. Section~\ref{s:SBR} studies the Sinai-Bowen-Ruelle measure, derives representation formulas for $S(\xi,\cdot)$, and reviews the theory of Bernoulli convolutions. Dually, Section~\ref{s:existence_localtime} investigates the occupation measures of $H$ and $\mathcal{T}$. We derive semi-explicit formulas for the Fourier transform of the occupation measure via Rademacher chaos expansions and combine them with Plancherel's theorem and the $L^2$-regularity theory of Bernoulli convolutions to establish the existence of square-integrable local times for drifted Takagi curves for almost every $\gamma\in(1/2,1)$. Finally, for the classical Takagi function, we use equidistribution modulo $\pi$, Koksma's theorem, and discrepancy estimates due to Erd\H{o}s and Tur\'an to prove the existence of square-integrable local times for $\gamma=2^{-1/m}$, $m\ge7$.

	\section{The curve as the attractor of a dynamical system}\label{s:attractor}
	
In this section, we introduce the dynamical systems framework that underlies the remainder of the paper. More precisely, we show that the graph of the Takagi-type curve can be interpreted as the global attractor of a dynamical system generated by the baker transform. This dynamical interpretation will later allow us to study both the associated Sinai-Bowen-Ruelle measure and the occupation measures of Takagi-type curves.
	
	Let $\gamma\in]\eh,1[$. In this paper, we investigate the fine geometric structure of the one-dimensional Takagi-type curve given by
	\usetagform{blue}
	\be
	\label{eq:DefinitionW}
	\cT(x)=\sum_{n=0}^{\infty}\gamma^n\Phi(2^n x),\qquad x\in[0,1],
	\ee
	where $\Phi(y)=d(y,\mathbb{Z})$, $y\in\mathbb{R}$. We begin by recalling that $\cT$ is H\"older continuous and determine its optimal H\"older exponent (see \cite{baranski15survey} for an overview).

	\begin{pr}\label{p:hoelderexponent}
		$\cT$ is H\"older continuous with exponent $-\frac{\log \gamma}{\log 2}$.
	\end{pr}
	
	\begin{proof}
		Let $x,y\in[0,1]$ and choose an integer $k\ge 0$ such that
		$$2^{-(k+1)}\leq |x-y|\leq 2^{-k}.$$
		Then we have, using the Lipschitz continuity of the distance function
		\bea
		|\cT(x) - \cT(y)| &\leq& \sum_{n=1}^k \gamma^n |d( 2^n x,\mathbb{Z}) - d(2^n y,\mathbb{Z})|
		+ 2\, \sum_{n=k+1}^\infty \gamma^n \\
		&\lesssim& \sum_{n=1}^k (2\gamma)^n |x-y| + \gamma^k \lesssim (2\gamma)^k\,\,2^{-k} + \gamma^k \simeq \gamma^k = 2^{-k \frac{\log \gamma}{\log \eh}}\\
		&\lesssim& |x-y|^{-\frac{\log \gamma}{\log 2}}.
		\eea
		This shows that $\frac{\log \gamma}{\log \eh}$ is an upper bound for the H\"older exponent of $\cT$. To see that it is also a lower bound, for $n\in\N$ choose $x_n = 0, y_n = 2^{-n}.$ Then we may write
		\bea
		|\cT(x_n) - \cT(y_n)| &=& \Big|\sum_{k=1}^\infty \gamma^k d(2^{k-n},\mathbb{Z})\Big|
		\\
		&=& \sum_{k=1}^{n-1} \gamma^k 2^{k-n} \simeq 2^{-n \frac{\log \gamma}{\log \eh}}
		=|x_n-y_n|^{-\frac{\log \gamma}{\log 2}}.
		\eea
		Since $|x_n-y_n|\to 0$ as $n\to\infty$, this shows that $-\frac{\log \gamma}{\log 2}$ is also a lower bound for the H\"older exponent of $\cT$. The argument can be extended to the other points in the interval.
	\end{proof}

	Our access to the analysis and geometry of $\cT$ is through the theory of dynamical systems. More precisely, we introduce a dynamical system on $[0,1]^2$, or equivalently on $\Om = \{0,1\}^\N\times\{0,1\}^\N$ whose global attractor is the graph of $\cT$. For $\omega\in\Omega$, we write for convenience $\om = ((\om_{-n})_{n\ge 0}, (\om_n)_{n\ge 1})$;so that $\Omega$ may be viewed as the space of $2$-dimensional sequences of Bernoulli random variables.  Denote by $\theta$ the canonical shift on $\Om$, given by
	$$
	\theta:\Om\to\Om,\quad \om\mapsto (\om_{n+1})_{n\in\Z}.
	$$
	We endow $\Om$ with the product $\sigma$-algebra, and the infinite product $\iota = \otimes_{n\in\Z} (\eh \delta_{\{0\}} + \eh \delta_{\{1\}})$under which the coordinate variables are independent Bernoulli random variables. We recall that $\theta$ is $\iota$-invariant.
	
	Now let
	$$
	T=(T_1,T_2):\Om \to [0,1]^2,\quad \om \mapsto (\sum_{n=0}^\infty \om_{-n} 2^{-(n+1)}, \sum_{n=1}^\infty \om_n 2^{-n}).
	$$


	Let $T_1$ and $T_2$ denote the first and second components of $T$, respectively.
	It is well known that $\iota$ is mapped by the transformation $T$ to $\lambda^2$ (i.e.~$\iota=\lambda^2\circ T$), the 2-dimensional Lebesgue measure. Moreover, the inverse of $T$,  which assigns to each point in $[0,1]^2$ its dyadic expansion, is uniquely defined except at dyadic points. For these we define the inverse to map to the sequences not converging to $0$. Let
	$$B =(B_1,B_2) = T\circ\theta \circ T^{-1}.$$
	We call $B=(B_1,B_2)$ the \emph{baker's transformation}. The $\theta$-invariance of $\iota$ directly translates into the $B$-invariance of $\lambda^2$:
	\begin{align}
		\label{eq:InvarianceBakerB}
		\lambda^2 \circ B^{-1} = (\lambda^2\circ T) \circ \theta^{-1}\circ T^{-1} = (\iota\circ \theta^{-1}) \circ T^{-1} = \iota \circ T^{-1} = \lambda^2.
	\end{align}
	For $(\xi, x)\in[0,1]^2$ let us note
	$$T^{-1}(\xi, x) = \big((\bxi_{-n})_{n\ge 0}, (\bx_n)_{n\ge 1}\big).$$
	Let us calculate the action of $B$ and its entire iterates on $[0,1]^2.$
	
	\begin{lem}\label{l:action_B}
		Let $(\xi, x) \in[0,1]^2$. Then for $k\ge 0$
		$$
		B^k(\xi, x)
		=
		\Big(2^k \xi\, (\textrm{mod } 1), \frac{\bxi_{-k+1}}{2} + \frac{\bxi_{-k+2}}{2^2}+\cdots+\frac{\bxi_0}{2^k} + \frac{x}{2^k}\Big),
		$$
		for $k\ge 1$
		$$ B^{-k}(\xi, x)
		=
		\Big(\frac{\xi}{2^k}+ \frac{\bx_1}{2^k} + \frac{\bx_2}{2^{k-1}}+\cdots+ \frac{\bx_k}{2}, 2^k x (\textrm{mod } 1)\Big).
		$$
	\end{lem}
	\pf
	By definition of $\theta^k$ for $k\ge 0$
	$$B^k(\xi,x) = \Big(\sum_{n\ge 0} \bxi_{-n+k} 2^{-(n+1)}, \frac{\bxi_{-k+1}}{2} + \frac{\bxi_{-k+2}}{2^2}+\cdots+\frac{\bxi_0}{2^k} + \sum_{n\ge 1} \bx_n 2^{-(k+n)}\Big).$$
	Now we can write
	$$\sum_{n\ge 0} \bxi_{-n+k} 2^{-(n+1)} = 2^k \xi (\mbox{mod } 1)
	\quad\textrm{and}\quad
	\sum_{n\ge 1} \bx_n 2^{-(k+n)} = \frac{x}{2^k}.$$
	This gives the first formula.
	For the second, note that by definition of $\theta^{-k}$ for $k\ge 1$
	$$B^{-k}(\xi, x) = \Big(\sum_{n\ge 0} \bxi_{-n} 2^{-(n+1+k)} + \frac{\bx_1}{2^k} + \frac{\bx_2}{2^{k-1}}+\cdots +\frac{\bx_k}{2}, \sum_{n\ge 1} \bx_{n+k} 2^{-n}\Big).$$
	Again, we identify
	$$\sum_{n\ge 1} \bx_{n+k} 2^{-n} = 2^k x (\textrm{mod } 1)
	\qquad\textrm{and}\qquad
	\sum_{n\ge 0} \bxi_{-n} 2^{-(n+1+k)} = \frac{\xi}{2^k}.$$
	\epf
	For $k\in\Z, (\xi, x)\in[0,1]^2$ we abbreviate the $k$-fold iterate of the baker transform of $(\xi,x)$ as
	$$B^{k}(\xi, x) = \big( B^{k}_1(\xi, x), B^{k}_2(\xi, x) \big) = (\xi_k, x_k),$$
	where for $k\ge 0$
	$$\xi_k = 2^k\xi (\textrm{mod } 1),
	\quad\textrm{and}\quad
	x_k = \frac{\bxi_{-k+1}}{2} + \frac{\bxi_{-k+2}}{2^2}+\cdots+\frac{\bxi_0}{2^k} + \frac{x}{2^k},$$
	and for $k\ge 1$
	$$\xi_{-k} = \frac{\xi}{2^k}+ \frac{\bx_1}{2^k} + \frac{\bx_2}{2^{k-1}}+\cdots + \frac{\bx_k}{2},
	\quad\textrm{and}\quad x_{-k} = 2^k x (\mbox{mod } 1).$$
	Following Baranski \cite{baranski02, baranski12,baranski14published}, Shen \cite{shen2017-Publication-of-2015}, Hunt \cite{hunt1998} and \cite{IdR2018}, we will next interpret the Takagi curve $\cT$ by a transformation on our base space $[0,1]^2$.
	Let
	\bea
	F: [0,1]^2\times\R &\to& [0,1]^2\times \R,\\
	(\xi,x, y)&\mapsto& \Big(B(\xi, x), \gamma y + \Phi(B_2(\xi, x))\Big).
	\eea
	Here we note $B = (B_1, B_2)$ for the two components of the baker transform $B$.
	
	For convenience, we extend $\cT$ from $[0,1]$ to $[0,1]^2$ by setting
	$$\cT(\xi, x) = \cT(x),\quad \xi, x\in[0,1].$$
	To see that the graph of $\cT$ is an attractor for $F$, the skew-product structure of $F$ with respect to $B$ plays a crucial role.
	\begin{lem}\label{l:attractor}
		For any $\xi, x\in[0,1]$ we have
		$$F\big(\xi, x, \cT(\xi, x)\big) = \big(B(\xi, x), \cT\big(B(\xi, x)\big)\big).$$
	\end{lem}
	
	\pf By the definition of the baker's transform we may write
	$$\cT(\xi, x) = \sum_{n=0}^\infty \gamma^n \Phi\big(B^{-n}_2(\xi,x)\big),\quad \xi, x\in[0,1].$$
	Hence, setting $k=n-1$, for $\xi, x\in[0,1]$
	\bea
	\cT\big(B_2(\xi,x)\big)
	&=& \sum_{n=0}^\un \gamma^n \Phi\big(B^{-n+1}_2(\xi,x)\big)\\
	&=& \Phi \big(B_2(\xi, x)\big) + \gamma \sum_{k=0}^\infty \gamma^k \Phi \big(B^{-k}_2(\xi,x)\big)\\
	&=& \Phi \big(B_2(\xi, x)\big) + \gamma \cT(x).
	\eea
	Hence by definition of $F$
	$$\big(B(\xi, x), \cT(B(\xi, x))\big)
	= \big(B(\xi, x), \cT(B_2(\xi, x))\big) = F\big(\xi, x, \cT(\xi, x)\big).$$
	\epf
	
	To assess the stability properties of the dynamical system generated by $F$, we calculate its Jacobian. We obtain for $\xi, x \in[0,1], y \in\R$
	\bea
	DF(\xi, x, y) &=&
	\left[
	\begin{array}{ccc}
		2&0&0\\
		0&\eh&0\\
		0& \eh \Phi' \big(B_2(\xi, x) \big)& \gamma
	\end{array}\right].
	\eea
	Hence the Lyapunov exponents of the dynamical system associated with $F$ are given by $2, \eh$, and $\gamma.$
	The corresponding invariant vector fields are given by
	$$\left(\begin{array}{c} 1\\0\\0\end{array}\right),\quad X(\xi,x)
	= \left(\begin{array}{c}0\\1\\ - \sum_{n=1}^\un \big(\frac{1}{2\gamma}\big)^n \Phi'\big(B^n_2(\xi,x)\big)\end{array}\right),\quad \left(\begin{array}{c} 0\\0\\1\end{array}\right),$$
	as is straightforwardly verified. Note that $X$ is well defined, since by our choice of $\gamma$ we have $2\gamma > 1$.
	Hence we have in particular for $\xi, x\in[0,1], y \in\R$
	$$DF(\xi, x, y) X(\xi, x) = \frac{1}{2}\,\, X\big(B(\xi, x)\big).$$
	Note that the vector $X$ spans an invariant stable manifold and does not depend on $y$.

	\section{The Sinai-Bowen-Ruelle measure}\label{s:SBR}
	
	Abbreviate $\kappa=\frac{1}{2\gamma}\in(0,1)$. Tsujii \cite{tsujii01} investigated the absolute continuity of the Sinai--Bowen--Ruelle (SBR) measure supported on the stable manifold
	\[
	S(\xi,x)=\sum_{n=1}^{\infty}\kappa^n\Phi'\big(B_2^n(\xi,x)\big),
	\qquad \xi,x\in[0,1],
	\]
	with respect to Lebesgue measure. His approach is based on the transversality of the family of functions
	\[
	x\longmapsto S(\xi,x)-S(\eta,x),\qquad \xi,\eta\in[0,1].
	\]
	In this section, we exploit the simple structure of the Takagi-type function $\cT$ to derive an explicit representation of $S$ in terms of a symmetric Bernoulli convolution, whose regularity properties have been studied extensively. This representation will play a key role throughout the paper.
	
	To describe the SBR measure associated with $F$, we first compute the action of $S$ on the  $\lambda^2$-measure preserving map $B$. For $\xi,x\in[0,1]$, we have

	\bea
	S(B(\xi, x)) &=& \sum_{n=1}^\un \kappa^n \Phi'\big(B^n_2\big(B_2(\xi, x)\big)\big)\\
	&=& \sum_{n=1}^\un \kappa^n \Phi'\big(B_2^{n+1}(\xi, x)\big)\\
	&=& \kappa^{-1} \sum_{k=1}^\un \kappa^k \Phi'\big(B^k_2(\xi, x)\big) - \Phi' \big(B_2(\xi, x)\big)\\
	&=& 2 \gamma S(\xi, x) - \Phi'\big(B_2(\xi, x)\big).
	\eea
	So we may define the Anosov skew product
	\begin{align*}
		\Gamma:[0,1]^2\times \R
		& \to[0,1]^2\times \R,
		\\
		(\xi, x, v)
		&\mapsto
		\big(B(\xi, x), 2 \gamma v - \Phi'\big(B_2(\xi, x)\big)\big).
	\end{align*}
	Then the equation just obtained yields the following result ({compare with Lemma \ref{l:attractor}}).
	\begin{lem}\label{l:SBR}
		For $\xi, x\in[0,1]$ we have
		$$\Gamma\big(\xi, x, S(\xi, x)\big) = \big(B(\xi, x), S(B(\xi, x))\big).$$
		The push-forward of the Lebesgue measure in $\R^2$ to the graph of $S$ given by
		$$\psi = \lambda^2\circ (\mbox{id}, S)^{-1}$$
		on $\cB([0,1]^2)\otimes \cB(\R)$ is $\Gamma$-invariant.
	\end{lem}
	
	\pf The first equation has been verified above. The $\Gamma$-invariance of $\psi$ is a direct consequence of the $B$-invariance of $\lambda^2.$ \epf
	
	Define $\pi_2:[0,1]^2\to[0,1], (\xi, x)\mapsto x$ and define the measure
	\begin{align}
		\label{eq:SBRforGamma}
		\mu & = \lambda^2\circ (\pi_2, S)^{-1}
	\end{align}
	on $\cB([0,1]^2)$. The measure $\mu$ is called the \emph{Sinai-Bowen-Ruelle measure} of $\Gamma$. Its marginals in $x\in[0,1]$ are denoted $\mu_x=\lambda \circ S(\cdot,x)^{-1}$.

	
	We now define a map on our probability space that exhibits certain increments of $S$ in a self similar way.
	Let
	
	$$G(\xi,x) = \sum_{n\in\Z} \kappa^{-n} \big[\Phi'\big(B^{-n}_2(\xi, x)\big) - \Phi' \big(B^{-n}_2(0, x)\big)\big],\quad \xi, x\in[0,1].$$
	
	Then we have the following simple relationship between $G$ and $S$.
	
	\begin{lem}\label{l:doubly_infinite_series_G}
		For $x, \xi, \eta\in [0,1]$ we have
		$$G(\xi, x) - G(\eta, x) = S(\xi,x) - S(\eta, x).$$
	\end{lem}
	
	\pf
	For $x, \xi, \eta\in[0,1]$, the equation
	\bea
	G(\xi, x) - G(\eta, x) &=& \sum_{n\in\Z} \kappa^{-n} \big[\Phi' \big(B^{-n}_2(\xi, x)\big) - \Phi' \big(B^{-n}_2(\eta, x)\big)\big]
	\\
	&=& \sum_{k=1}^\infty \kappa^k \big[\Phi' \big(B^{k}_2(\xi, x)\big) - \Phi' \big(B^{k}_2(\eta, x)\big)\big]\\
	&=& S(\xi, x) - S(\eta, x),
	\eea
	holds. Here we used that the first sum for non-negative integers $n$ is zero.
	
	This completes the proof. \epf
	
	The following result describes the scaling properties of $G$.
	
	\begin{lem}
		\label{l:scaling_G}
		For $\xi, x\in[0,1]$ we have
		$$G\big(B^{-1}(\xi, x)\big) = \kappa G(\xi,x).$$
	\end{lem}
	
	\pf
	Note that by definition, setting $n+1 = k$, for $\xi, x\in[0,1]$
	\bea
	G\big(B^{-1}(\xi, x)\big)
	&=&\sum_{n\in\Z} \kappa^{-n} \big[\Phi' \big( B^{-n-1}(\xi, x)\big) - \Phi' \big( B^{-n}(B^{-1}_1(\xi,x), 0)\big)\big]\\
	&&+ \sum_{n\in\Z} \kappa^{-n} \big[\Phi' \big( B^{-n-1}(\xi, x)\big) - \Phi' \big( B^{-n-1}(\xi, 0)\big)\big]\\
	&&+ \sum_{n\in\Z} \kappa^{-n} [\Phi' \big( B_2^{-n-1}(\xi, 0)\big) - \Phi' \big( B^{-n}_2(B_1^{-1}(\xi,x), 0)\big)\big] \\
	&=& \kappa G(\xi,x) + \sum_{k=1}^\infty \kappa^{k} [\Phi'(B_2^{k-1}(\xi,0))-\Phi'(B_2^k(B_1^{-1}(\xi,x),0))]\\
	&=& \kappa G(\xi,x).
	\eea
	For the last equality, note that $\Phi'$ is constant on the two halves of the unit interval, and that
	$$B_2^{k-1}(\xi,0) = \frac{\bxi_0}{2^{k-1}}+\cdots + \frac{\bxi_{-k+2}}{2}$$
	and
	$$B_2^k(B_1^{-1}(\xi,x),0) = \frac{\bx_1}{2^k} + \frac{\bxi_0}{2^{k-1}}+\cdots+\frac{\bxi_{-k+2}}{2}$$
	belong to the same half.
	This provides the claimed identity.\epf
	
	We conclude this section by deriving a representation of $S$ that forms the starting point for our analysis of its transversality properties. To this end, fix $\xi,\eta\in[0,1]$. We recursively define the sequence of disagreement times of the dyadic components of $\xi$ and $\eta$. For $n\in\mathbb N$, let	
	\begin{align}
		\tau_{1}=\inf\{\ell\geq 0:\bar{\xi}_{-\ell}\not= \bar{\eta}_{-\ell}\},
		\text{ and }
		\tau_{n+1}=\inf\{\ell>\tau_n:\bar{\xi}_{-\ell}\not= \bar{\eta}_{-\ell}\},
	\end{align}
	
	and for $x\in[0,1]$
	\begin{align*}\label{defg1}
		g(x):=&\sum^{\infty}_{m=0}\kappa^m\big[
		\Phi' \big(B_2^{m}(0,\frac{1+x}{2})\big)-\Phi'\big(B_2^{m}(0,\frac{x}{2})\big)\big]\\
		=&\sum^{\infty}_{m=0}\kappa^m\big[
		\Phi' \big(\frac{1+x}{2^{m+1}}\big)-\Phi' \big(\frac{x}{2^{m+1}}\big)\big].
	\end{align*}
	
	We have the following result.
	
	\begin{pr}\label{p:representation_S}
		Let $\xi,\eta,x\in[0,1]$. Then
		\be\label{e:representation_G}
		S(\xi,x)-S(\eta,x)=\sum^{\infty}_{\ell=1}\kappa^{\tau_\ell+1}\,(-1)^{(1-\bar{\xi}_{-\tau_{\ell}})} g\big(B_2^{\tau_\ell}(\xi,x)\big).
		\ee
	\end{pr}
	
	\pf
	It follows from the definition of $\tau_n, n\in\N,$ that $\xi$ can be written
	$$
	\xi=(\beeta_0,\ldots,\beeta_{-\tau_1+1},\bxi_{-\tau_1},\beeta_{-\tau_1-1},\cdots,\beeta_{-\tau_2+1},\bxi_{-\tau_2},\beeta_{-\tau_2-1},\cdots\cdots,\beeta_{-\tau_n+1},\bxi_{-\tau_n},\beeta_{-\tau_n-1},\cdots).
	$$
	For $n\in \N$ let $\xi^n$ be the sequence which up to $\tau_n$ represents the dyadic expansion of $\xi$, and then switches to the representing sequence of $\eta.$ Then for $n\in\N$ we have
	$$
	\xi^n=(\beeta_0,\ldots,\beeta_{-\tau_1+1},\bxi_{-\tau_1},\beeta_{-\tau_1-1},\cdots,\beeta_{-\tau_2+1},\bxi_{-\tau_2},\beeta_{-\tau_2-1},\cdots\cdots,\beeta_{-\tau_n+1},\bxi_{-\tau_n},\beeta_{-\tau_n-1},\cdots,\beeta_{-m},\cdots).
	$$
	Note that $\lim_{n\to\infty} \xi^n = \xi.$ We can therefore rewrite the left hand side of \eqref{e:representation_G} in the telescoping form
	$$
	S(\xi,x)-S(\eta,x)=\sum^{\infty}_{l=1}\Big(S(\xi^\ell,x)-S(\xi^{\ell-1},x)\Big),
	$$
	where $\xi^0=\eta.$

	For $\ell\in \N$ let us calculate $S(\xi^\ell,\cdot)-S(\xi^{\ell-1},\cdot)$. Since $\bxi^{\ell}_{-k}=\bxi^{\ell-1}_{-k}$ for $k\leq \tau_\ell-1$, we have
	\begin{align}
		S(\xi^\ell,x)-S(\xi^{\ell-1},x)
		=&\sum^{\infty}_{n=1}\kappa^n\Big[
		\Phi' \big(B_2^n(\xi^{\ell},x)\big)-\Phi' \big(B_2^n(\xi^{\ell-1},x)\big)\Big]\notag
		\\
		=&\sum^{\infty}_{n=\tau_\ell +1}\kappa^n\Big[
		\Phi' \big(B_2^n(\xi^{\ell},x)\big)-\Phi' \big(B_2^n(\xi^{\ell-1},x)\big)\Big]\notag\\
		=&\kappa^{\tau_\ell}\sum^{\infty}_{m=1}\kappa^m\Big[
		\Phi' \Big(B_2^m\big(B^{\tau_\ell}(\xi^{\ell},x)\big)\Big)
		-\Phi' \Big(B_2^m\big(B^{\tau_\ell}(\xi^{\ell-1},x)\big)\Big)\Big].\notag
	\end{align}
	Now
	$$
	B_2^{\tau_\ell}(\xi^{\ell},x)=B_2^{\tau_\ell}(\xi,x)=B_2^{\tau_\ell}(\xi^{\ell-1},x).
	$$
	In case $\bxi_{-\tau_{\ell}}=1$ we have
	$$
	B_2^1\big(B^{\tau_\ell}(\xi^{\ell},x)\big)
	=B_2^{\tau_\ell+1}(\xi^{\ell},x)
	=\frac{1+B_2^{\tau_\ell}(\xi^{\ell},x)}{2}
	=\frac{1+B_2^{\tau_\ell}(\xi,x)}{2},
	$$
	while
	$$
	B_2^1\textbf{}\big(B^{\tau_\ell}(\xi^{\ell-1},x)\textbf{}\big)
	=B_2^{\tau_\ell+1}(\xi^{\ell-1},x)
	=\frac{B_2^{\tau_\ell}(\xi,x)}{2}.
	$$
	In case $\bxi_{-\tau_{\ell}}=0$, we have in contrast
	$$
	B_2^1\big(B^{\tau_\ell}(\xi^{\ell},x)\big)
	=B_2^{\tau_\ell+1}(\xi^{\ell},x)
	=\frac{B_2^{\tau_\ell}(\xi^{\ell},x)}{2}
	=\frac{B_2^{\tau_\ell}(\xi,x)}{2},
	$$
	while
	$$
	B_2^1\textbf{}\big(B^{\tau_\ell}(\xi^{\ell-1},x)\textbf{}\big)
	=B_2^{\tau_\ell+1}(\xi^{\ell-1},x)
	=\frac{1+B_2^{\tau_\ell}(\xi^{\ell},x)}{2}=\frac{1+B_2^{\tau_\ell}(\xi,x)}{2}.
	$$
	Hence we may write by definition of $g$
	\begin{align*}
		S(\xi^\ell,x)-&S(\xi^{\ell-1},x)
		\\
		=
		&\,(-1)^{(1-\bxi_{-\tau_{\ell}})}\,\kappa^{\tau_\ell} \sum^{\infty}_{m=1} \kappa^m
		\Big[
		\Phi' \Big(B_2^{m-1}\Big(\frac{1+B_2^{\tau_\ell}(\xi,x)}{2}\Big)\Big)
		-\Phi' \Big(B_2^{m-1}\Big(\frac{B_2^{\tau_\ell}(\xi,x)}{2}\Big)\Big)\Big]\notag
		\\
		=& \kappa^{\tau_\ell+1}\,(-1)^{(1-\bxi_{-\tau_{\ell}})}\,g\big(B_2^{\tau_\ell}(\xi,x)\big).
	\end{align*}
	
	Hence we obtain the claimed representation
	$$
	S(\xi,x)-S(\eta,x)=\sum^{\infty}_{\ell=1}\kappa^{\tau_\ell+1}\,(-1)^{(1-\bxi_{-\tau_{\ell}})}\, g\big(B_2^{\tau_\ell}(\xi,x)\big), \,\,\,\xi,\eta,x\in [0,1].
	$$\epf
	
	Let us calculate $g$. Here, as opposed to the trigonometric case in the stable manifold of Weierstrass curves (see \cite{imkellerdosreispamenreveillac2020}), the simplicity of $\phi$ implies the following surprising identity.
	
	\begin{lem}\label{l:structure_g}
		We have $g(x) = -2, x\in[0,1].$
	\end{lem}
	
	\pf Let us inspect the first term in the series decomposition of $g$. Here we have $B_2^0(0,\frac{1+x}{2}) = \frac{1+x}{2} \in[\eh,1],$
	while $B_2^0(0,\frac{x}{2}) = \frac{x}{2} \in[0,\eh].$ Hence the contribution of the first term is $-1-1=-2$. For $m\ge 1$ we have in contrast that both
	$B_2^m(0,\frac{1+x}{2}) = \frac{1+x}{2^{m+1}}$ and $B_2^m(0,\frac{x}{2}) = \frac{x}{2^{m+1}}$ belong to $[0,\eh].$ Hence the contribution of terms of order $m\ge 1$ vanishes. This implies the claimed identity.\epf
	
	Lemma \ref{l:structure_g} gives the following simplification of the representation formula of Proposition \ref{p:representation_S}.
	
	\begin{co}\label{c:simple_representation_S}
		Let $\xi,\eta,x\in[0,1]$. Then
		\be\label{e:representation2_G}
		S(\xi,x)-S(\eta,x)=-2\,\sum^{\infty}_{\ell=1}\kappa^{\tau_\ell+1}\,(-1)^{(1-\bar{\xi}_{-\tau_{\ell}})}.
		\ee
	\end{co}
	
	\pf This follows by combining Lemma \ref{l:structure_g} and Proposition \ref{p:representation_S}. \epf
	
	Let us finally extend this property to $S(\xi,.)$ for $\xi\in[0,1].$ First observe that for any $x\in[0,1]$ we have
	\be\label{e:representation_S_at0}
	S(0,x) = \sum_{n=1}^\infty \kappa^n \Phi'(\frac{x}{2^n}) = \sum_{n=1}^\infty \kappa^n = \frac{\kappa}{1-\kappa}.
	\ee
	Now denote by $\tau^0_n, n\in\mathbb{N}$, the sequence of stopping times described above for the particular case $\eta=0$. Then we obtain
	
	\begin{co}\label{c:simple_representation_S_oneterm}
		Let $\xi, x\in[0,1]$. Then
		\be\label{e:representation_S_oneterm}
		S(\xi,x) = \frac{\kappa}{1-\kappa} - 2\,\sum^{\infty}_{\ell=1}\kappa^{\tau^0_\ell+1} = \sum_{n=0}^\infty \big( 1_{\{\bxi_{-n}=0\}}-1_{\{\bxi_{-n}=1\}}\big) \kappa^{n+1}.
		\ee
	\end{co}
	
	\pf Combine Corollary \ref{c:simple_representation_S} with \eqref{e:representation_S_at0}, and note that by definition $\bxi_{-\tau^0_k}=1$ for all $k\ge 1.$ \epf

	The final identity in the preceding corollary shows that $S$ can be identified with a symmetric Bernoulli convolution in the sense of Erd\H{o}s \cite{erdoes39}.
	
	In Corollary~\ref{c:simple_representation_S_oneterm}, we identified $S$ with a symmetric Bernoulli convolution. Bernoulli convolutions have occupied a central place in rough geometry, harmonic analysis, and dynamical systems since the pioneering works of Wintner \cite{wintner35}, Kershner and Wintner \cite{kershnerwintner35}, and Erd\H{o}s \cite{erdoes39} in the 1930s. They have subsequently been studied by many authors, including Brauer \cite{brauer51}, Solomyak \cite{solomyak95}, and Shmerkin \cite{shmerkin14}, and have attracted renewed attention over the past three decades. We refer to the surveys of Peres \cite{peres98} and Varj\'u \cite{varju18} for comprehensive accounts.
	
	A central problem in the theory is to determine for which values of $\kappa\in(0,1)$ the random series
	\[
	Z(\kappa)=\sum_{n=1}^{\infty}X_n\kappa^n,
	\]
	where $(X_n)_{n\in\mathbb N}$ is a sequence of i.i.d.\ symmetric Bernoulli random variables, has a distribution that is absolutely continuous with respect to Lebesgue measure. From the earliest works onwards, this question has been approached through Fourier analysis. The characteristic function (Fourier transform) of $Z(\kappa)$ is given by
	\[
	\phi_Z(u)=\prod_{n=1}^{\infty}\cos(u\kappa^n),
	\qquad u\in\mathbb R,
	\]
	and Plancherel's theorem implies that the distribution of $Z(\kappa)$ possesses a square-integrable density whenever $\phi_Z\in L^2(\mathbb R)$. One of the features of Bernoulli convolutions is that this property depends in a non-trivial way on the algebraic property of the parameter $\kappa$. Early discoveries about this problem lead to the conclusion that either the law is absolutely continuous or singular with respect to Lebesgue measure.  
	
	The first positive result is due to Wintner \cite{wintner35}, who proved absolute continuity for $\kappa=2^{-1/m},\,\,\, m\in\mathbb N,$ a parameter regime that will reappear in our study of local times. Shortly afterwards, Erd\H{o}s \cite{erdoes39} showed that the distribution is singular whenever $\kappa^{-1}$ is a Pisot number  (i.e. its Galois conjugates, the roots of its minimal polynomial over the complex field are all contained inside the unit disc.), thereby revealing the unexpectedly strong influence of algebraic properties. A major breakthrough came later, when Solomyak \cite{solomyak95} proved that Bernoulli convolutions are absolutely continuous for almost every $\kappa\in(\frac12,1)$. Despite this result, $\frac12$ remains a cluster point of parameters for which the corresponding Bernoulli convolution is singular. We refer the reader to the survey by Varj\'u \cite{varju18} for a comprehensive introduction to the subject and an account of many of its major developments. In the present paper, we exploit the representation established in Corollary~\ref{c:simple_representation_S_oneterm} to transfer the above results on Bernoulli convolutions to corresponding properties of occupation measures and local times of the Takagi curve $\cT$.

	\section{The existence of a local time for $\cT$}\label{s:existence_localtime}
	In this section, we establish the existence of a square-integrable local time for the occupation measure associated with $\cT$. As in the previous section, the proof relies on a Fourier-analytic criterion due to Berman.

	\subsection{Occupation measures and local times}
In this subsection, we recall the notions of occupation measures and local times, together with Berman's Fourier criterion, which forms the analytical foundation of our approach.
	
Let $f:[0,1]\rightarrow\mathbb{R}$ be a measurable function. The
		\textsl{occupation measure} of $f$ is the probability measure
		$\mu_f$ defined by
		\[
		\mu_f(A)
		=
		\lambda\bigl(\{x\in[0,1]:f(x)\in A\}\bigr),
		\quad
		A\in\mathcal{B}(\mathbb{R}),
		\]
		where $\lambda$ denotes the Lebesgue measure on $[0,1]$. Equivalently,
		its distribution function is
		$$
		F_f(y)
		=
		\mu_f((-\infty,y])
		=
		\int_0^1
		\mathbf{1}_{\{f(x)\le y\}}
		\,\mathrm{d}x.
		$$
	The occupation measure is characterized by the identity
		\begin{equation}\label{eq:occupation_formula}
			\int_0^1 g(f(x))\,\mathrm{d}x
			=
			\int_{\mathbb{R}} g(y)\,\mu_f(\mathrm{d}y),
		\end{equation}
		which holds for every bounded Borel measurable function $g$.
		
		Taking $g(y)=e^{iuy}$ yields the Fourier transform (or characteristic
		function) of $\mu_H$,
		\begin{equation}\label{eq:FToccupation}
			\widehat{\mu}_f(u)
			=
			\int_{\mathbb{R}}e^{iuy}\,\mu_f(\mathrm{d}y)
			=
			\int_0^1e^{iuf(x)}\,\mathrm{d}x,
			\quad u\in\mathbb{R}.
	\end{equation}
	If $\mu_f$ is absolutely continuous with respect to Lebesgue measure,
		its Radon-Nikodym derivative
		\[
		L_f(y)
		=
		\frac{\mathrm{d}\mu_f}{\mathrm{d}y}(y)
		\]
		is called the \textsl{local time} of $f$. The following classical criterion provides a sufficient condition for the existence of square-integrable
		local times.
		
		\begin{lem}[Berman {\cite{berman69}}]
			\label{lem:berman}
			Let $f:[0,1]\to\mathbb{R}$ be measurable and let $\mu_f$ denote its
			occupation measure. Then:
			\begin{enumerate}
				\item[(i)] If
				\[
				\int_{\mathbb R}
				|\widehat{\mu_f}(u)|^2\,\mathrm{d}u<\infty,
				\]
				then $\mu_f$ is absolutely continuous with respect to Lebesgue measure,
				and its local time belongs to $L^2(\mathbb R)$.
				\item[(ii)] If
				\[
				\int_{\mathbb R}
				|\widehat{\mu_f}(u)|\,\mathrm{d}u<\infty,
				\]
				then the local time admits a continuous version.
			\end{enumerate}
		\end{lem}
		Lemma \ref{lem:berman} reduces the existence of local times to obtaining suitable estimates for the Fourier transform of the occupation measure. We now derive such estimates for the Takagi curve.

	\subsection{Representation and scaling properties of $H$}
	Our strategy is to estimate the Fourier transform of the occupation measure through an auxiliary function that admits a convenient representation in terms of the Bernoulli convolution $S$. To this end, we introduce the function
	
	\begin{align}
		\label{eq:FunctionH}
		H(\xi,x)
		=
		\sum_{n\in\mathbb Z}
		\gamma^n
		\Big[
		\Phi\big(B_2^{-n}(\xi,x)\big)
		-
		\Phi\big(B_2^{-n}(\xi,0)\big)
		\Big],
		\qquad
		\xi,x\in[0,1],
	\end{align}
	The first step is to relate $H$ to the Takagi curve and its stable manifold. We shall prove that
	\[
	H(\xi,x)
	=
	\cT(x)-\int_0^xS(\xi,z)\,\mathrm dz
	=
	\cT(x)-xS(\xi,0),
	\]
	since $S(\xi,\cdot)$ is constant. We then show that $H$ admits the representation
	\[
	H(\xi,y)
	=
	\kappa
	\sum_{\ell=1}^\infty
	\gamma^{\sigma_\ell^0}
	S(B_1^{-\sigma_\ell^0+1}(\xi,y),0),
	\]
	which expresses the geometry of the Takagi curve entirely in terms of Bernoulli convolutions $S$. Thus $H$ plays for occupation measures the same role that $G$ plays for the SBR measure in the previous section, and this representation forms the basis of the subsequent Fourier analysis.

	\begin{lem}\label{l:doubly_infinite_series}
		For $x, y, \xi\in [0,1]$ we have
		$$H(\xi, y) - H(\xi, x) = \cT(y) - \cT(x) - \int_x^y S(\xi, z) \dz.$$
		For $x, \xi,\eta \in [0,1]$ we have
		$$H(\eta, x) - H(\xi, x) = \int_0^x \big(S(\xi, z)-S(\eta,z)\big) \dz= x (S(\xi,0)-S(\eta,0)).$$
	\end{lem}
	
	\pf
	For $x, y, \xi\in[0,1]$ we have indeed (recall $\kappa=\frac1{2\gamma}$)
	\bea
	H(\xi, y) - H(\xi, x)
	&=& \sum_{n\in\Z} \gamma^n \big[\Phi\big( B^{-n}_2(\xi, y) \big)- \Phi\big( B^{-n}_2(\xi, x)\big)\big]
	\\
	&=& \sum_{n=0}^\infty \gamma^n \big[\Phi\big( B^{-n}_2(\xi, y)\big) - \Phi \big( B^{-n}_2(\xi, x)\big)\big]
	\\
	&&\hspace{1.5cm} + \sum_{k=1}^\infty \gamma^{-k} \big[\Phi \big( B^{k}_2(\xi, y)\big) - \Phi \big( B^{k}_2(\xi, x)\big)]
	\\
	&=& \cT(y)-\cT(x) - \int_x^y\, \sum_{k=1}^\infty (2 \gamma)^{-k} \Phi' \big( B^k_2(\xi,z)\big) \dz\\
	&=& \cT(y)-\cT(x) - \int_x^y  \sum_{k=1}^\infty \kappa^k \Phi' \big(B^k_2(\xi,z)\big) \dz\\
	&=& \cT(y) - \cT(x) - \int_x^y S(\xi, z) \dz.
	\eea				
	To argue for the second equation, note that for $x, \xi,\eta\in[0,1]$ we have
	\bea
	H(\eta, x) - H(\xi, x)
	&=& \sum_{n\in\Z} \gamma^n \Big\{\big[\Phi\big( B^{-n}_2(\eta, x) \big)-\Phi \big( B^{-n}_2(\eta, 0)\big)\big]
	\\
	&&\hspace{2cm} - \big[\Phi\big( B^{-n}_2(\xi, x)\big) - \Phi \big( B^{-n}_2(\xi, 0)\big)\big]\Big\}
	\\
	&=& \sum_{k=1}^\infty \gamma^{-k} \Big\{\big[\Phi \big( B^{k}_2(\eta, x)\big) -\Phi \big( B^{k}_2(\eta, 0)\big)\big]
	\\
	&&\hspace{2cm} - \big[\Phi\big( B^{k}_2(\xi, x)\big) - \Phi \big( B^{k}_2(\xi, 0)\big)\big]\Big\}
	\\
	&=& -\int_0^x  \sum_{k=1}^\infty (2 \gamma)^{-k} \big[\Phi' \big( B^k_2(\eta,z)\big)-\Phi' \big( B^k_2(\xi,z)\big)\big] \dz\\
	&=& -\int_0^x  \sum_{k=1}^\infty \kappa^k \big[\Phi' \big( B^k_2(\eta,z)\big)-\Phi' \big( B^k_2(\xi,z)\big)\big] \dz\\
	&=&  \int_0^x [S(\xi, z)-S(\eta,z)] \dz.
	\eea
	This completes the proof of the second equation. \epf
	
	We next address the scaling properties of $H$.
	\begin{lem}\label{l:scaling_H}
		For $\xi, x\in[0,1]$ we have 
		$$H\big( B(\xi, x) \big) = \gamma H(\xi,x) + \frac{\bxi_0}{2} S(2\xi,0).$$
		Consequently, for $\xi,x,y\in[0,1]$
		$$H(B(\xi,y))-H(B(\xi,x)) = \gamma[H(\xi,y)-H(\xi,x)].$$
	\end{lem}
	
	\pf
	Let $\xi, x\in[0,1]$. By definition and setting $n-1 = k$ we obtain, defining $\hat{\xi}$ to be represented by the dyadic sequence $(0,\bxi_{-1},\bxi_{-2},\cdots)$,
	\begin{align*}
		H(B(\xi, x)) =& \sum_{n\in\Z} \gamma^n \big[\Phi \big( B_{2}^{-n+1}(\xi, x)\big) - \Phi \big( B^{-n}_{2}(B_1(\xi,x), 0)\big)\big]
		\\
		=& \sum_{n\in\Z} \gamma^n \big[\Phi \big( B_{2}^{-n+1}(\xi, x)\big) - \Phi \big( B_{2}^{-n+1}(\xi, 0)\big)\big]\\
		&+ \sum_{n\in\Z} \gamma^n \big[\Phi \big( B_{2}^{-n+1}(\xi, 0)\big) - \Phi \big( B_{2}^{-n}(B_1(\xi,x), 0)\big)\big]\\
		=& \gamma \sum_{k\in\Z} \gamma^k \big[\Phi \big( B^{-k}(\xi, x)\big) - \Phi \big( B^{-k}(\xi, 0)\big)\big] \\
		&+ \sum_{k=1}^\infty \gamma^{-k} [\Phi \big( B_2^{k+1}(\xi, 0)\big) - \Phi \big( B^{k+1}_2(\hat{\xi}, 0)\big)\big] \\
		=& \gamma H(\xi,x) + \bxi_0 \sum_{k=1}^\infty \gamma^{-k} 2^{-k-1} \Phi'(B_2^{k+1}(\hat{\xi},0))\\
		=& \gamma H(\xi,x) + \frac{\bxi_0}{2} \sum_{k=1}^\infty \kappa^k \Phi'(B_2^k(2\xi,0))\\
		=& \gamma H(\xi,x)+ \frac{\bxi_0}{2} S(2\xi,0).
	\end{align*}
	The last equation follows from the fact that the second term in the first formula only depends on $\xi.$
	\epf

	We conclude this subsection by deriving a representation of the increments of $H$ that complements Proposition \ref{p:representation_S}. Whereas Proposition~\ref{p:representation_S} expresses increments of the stable manifold $S$ in terms of the function $g$, the representation below describes increments of $H$ in terms of the Bernoulli convolution structure underlying the Takagi curve. This identity provides the key ingredient in our study of occupation measures and the existence of square-integrable local times. To this end, fix $x,y\in[0,1]$. We recursively define the sequence of disagreement times of the dyadic components of $x$ and $y$. For $n\in\mathbb N$, let

	\begin{align}
		\sigma_{1}=\inf\{\ell\geq 1:\bar{x}_{\ell}\not= \bar{y}_{\ell}\},
		\text{ and }
		\sigma_{n+1}=\inf\{\ell>\sigma_n:\bar{x}_{\ell}\not= \bar{y}_{\ell}\}.
	\end{align}
	In case $x=0$ we denote this sequence by $(\sigma_\ell^0)_{\ell\in\N}$.
	
	
	In analogy to the increment functions $g$ in Section \ref{s:SBR} we define the following function
	$$ h(\xi,x) := \sum_{m=0}^\infty \gamma^{-m-1} \big[\Phi(B_2^m(\xi,\frac{1+x}{2}))-\Phi(B_2^m(\xi,\frac{x}{2}))\big],\quad \xi, x\in[0,1].
	$$
	
	We have the following result.
	
	\begin{pr}\label{p:representation_H}
		Let $\xi, x, y\in[0,1]$. Then
		\be\label{e:representation_H}
		H(\xi,y)-H(\xi,x)=\sum^{\infty}_{\ell=1}\gamma^{\sigma_\ell}\,(-1)^{(1-\bar{y}_{\sigma_{\ell}})} h(B_1^{-\sigma_\ell+1}(\xi,y),B_2^{-\sigma_\ell}(\xi,x)).
		\ee
	\end{pr}
	
	\pf
	It follows from the definition of $\sigma_n, n\in\N,$ that $y$ can be written
	$$
	y=(\bx_1,\ldots,\bx_{\sigma_1-1},\by_{\sigma_1},\bx_{\sigma_1+1},\cdots,\bx_{\sigma_2-1},\by_{\sigma_2},\bx_{\sigma_2+1},\cdots\cdots,\bx_{\sigma_n-1},\by_{\sigma_n},\bx_{\sigma_n+1},\cdots).
	$$
	For $n\in \N$ let $y^n$ be the sequence which up to $\sigma_n$ represents the dyadic expansion of $y$, and then switches to the representing sequence of $x.$ Then for $n\in\N$ we have
	$$
	y^n=(\bx_1,\ldots,\bx_{\sigma_1-1},\by_{\sigma_1},\bx_{\sigma_1+1},\cdots,\bx_{\sigma_2-1},\by_{\sigma_2},\bx_{\sigma_2+1},\cdots\cdots,\bx_{\sigma_n-1},\by_{\sigma_n},\bx_{\sigma_n+1},\cdots,\bx_{m},\cdots).
	$$
	Note that $\lim_{n\to\infty} y^n = y.$ Hence the left hand side of \eqref{e:representation_H} has the telescoping representation
	$$
	H(\xi,y)-H(\xi,x)=\sum^{\infty}_{l=1}\Big(H(\xi,y^\ell)-H(\xi,y^{\ell-1})\Big),
	$$
	where $y^0=x.$

	For $\ell\in \N$ let us calculate $H(\cdot,y^\ell)-H(\cdot, y^{\ell-1})$. Since $B_2^{-n}(\xi,y^\ell) = B_2^{-n}(\xi,y^{\ell-1})$ for $n\ge \sigma_\ell$ we have
	\begin{align}
		H(\xi,y^\ell)-H(\xi,y^{\ell-1})
		=&\sum_{n\in\Z}\gamma^n\Big[
		\Phi \big(B_2^{-n}(\xi,y^{\ell})\big)-\Phi \big(B_2^{-n}(\xi,y^{\ell-1})\big)\Big]\notag
		\\
		=&\sum_{n<\sigma_\ell}\gamma^n\Big[
		\Phi \big(B_2^{-n}(\xi,y^{\ell})\big)-\Phi \big(B_2^{-n}(\xi,y^{\ell-1})\big)\Big]\notag\\
		=&\gamma^{\sigma_{\ell}}\sum^{\infty}_{m=1}\gamma^{-m}\Big[
		\Phi \Big(B_2^m\big(B^{-\sigma_\ell}(\xi,y^{\ell})\big)\Big)
		-\Phi \Big(B_2^m\big(B^{-\sigma_\ell}(\xi,y^{\ell-1})\big)\Big)\Big]\notag\\
		=&\gamma^{\sigma_\ell} \sum^{\infty}_{k=0}\gamma^{-k-1}\Big[
		\Phi \Big(B_2^k\big(B^1(B^{-\sigma_\ell}(\xi,y^{\ell}))\big)\Big)
		-\Phi \Big(B_2^k\big(B^1(B^{-\sigma_\ell}(\xi,y^{\ell-1}))\big)\Big)\Big]\notag.
	\end{align}
	
	Now note that by definition $B^{-\sigma_\ell}(\xi,y^\ell)$ corresponds to the doubly infinite sequence
	$$(\cdots,\bxi_0,\by_1,\cdots,\by_{\sigma_\ell-1},\by_{\sigma_\ell},\bx_{\sigma_\ell+1},\cdots),$$
	while $B^{-\sigma_\ell}(\xi,y^{\ell-1})$ is represented by
	$$(\cdots,\bxi_0,\by_1,\cdots,\by_{\sigma_\ell-1},\bx_{\sigma_\ell},\bx_{\sigma_\ell+1},\cdots),$$
	where $\by_{\sigma_\ell}$ resp. $\bx_{\sigma_\ell}$ take position $0$ in the sequences. Hence in case $\by_{\sigma_\ell}=1$ we may continue to write
	\bea
	&& H(\xi,y^\ell)-H(\xi,y^{\ell-1}) \\
	&&\hspace{1cm} \gamma^{\sigma_\ell}\,(-1)^{(1-\by_{\sigma_\ell})}\,h(B_1^{-\sigma_\ell+1}(\xi,y),B_2^{-\sigma_\ell}(\xi,x)),
	\eea
	while in case $\by_{\sigma_\ell}=0$ the sign on the rhs just changes.
	This provides the claimed formula.

	\epf
	
	Let us next calculate $h$.
	
	\begin{lem}\label{l:structure_h}
		We have
		$$h(\xi,x) = 2\kappa[\Phi(\frac{1+x}{2})-\Phi(\frac{x}{2})] + \kappa S(\xi,0) = -2\kappa\,x + \kappa S(\xi,0),\quad  \xi,x\in[0,1].$$
	\end{lem}
	
	\pf
	We may write
	\bea
	h(\xi,x) &=& 2\kappa\,[\Phi(\frac{1+x}{2})-\Phi(\frac{x}{2})]\\
	&&+ \sum_{m=1}^\infty \gamma^{-m-1} \big[\Phi(B_2^m(\xi,\frac{1+x}{2}))-\Phi(B_2^m(\xi,\frac{x}{2}))\big]\\
	&=& 2\kappa\,[\Phi(\frac{1+x}{2})-\Phi(\frac{x}{2})]\\
	&&+ \sum_{m=1}^\infty \gamma^{-m-1} \big[\Phi(\frac{1+x}{2^{m+1}}+\frac{\bxi_0}{2^m}+\cdots+\frac{\bxi_{-m+1}}{2})-\Phi(\frac{x}{2^{m+1}}+\frac{\bxi_0}{2^m}+\cdots+\frac{\bxi_{-m+1}}{2})\big]\\
	&=&  2\kappa [\Phi(\frac{1+x}{2})-\Phi(\frac{x}{2})]\\
	&&+ \sum_{m=1}^\infty \gamma^{-m-1} 2^{-m-1}\, \Phi'(\frac{\bxi_0}{2^m}+\cdots+\frac{\bxi_{-m+1}}{2})\\
	&=& 2\kappa\,[\Phi(\frac{1+x}{2})-\Phi(\frac{x}{2})] + \kappa S(\xi,0).
	\eea
	Finally note that for $x\in[0,1]$ we have $\Phi(\frac{1+x}{2})-\Phi(\frac{x}{2}) = -x$. This completes the proof.\epf
	
	Using Lemma \ref{l:structure_h}, we may write the formula of Proposition \ref{p:representation_H} in a more explicit way.
	
	\begin{pr}\label{p:representation_H_explicit}
		Let $\xi, x, y\in[0,1]$. Then
		\bean\label{e:representation_H_explicit}
		H(\xi,y)-H(\xi,x)&=& - 2\kappa\,\sum^{\infty}_{\ell=1}\gamma^{\sigma_\ell}\,(-1)^{(1-\bar{y}_{\sigma_{\ell}})} B_2^{-\sigma_\ell}(\xi,x)\\
		\nonumber && + \kappa \sum^{\infty}_{\ell=1}\gamma^{\sigma_\ell}\,(-1)^{(1-\bar{y}_{\sigma_{\ell}})} S(B_1^{-\sigma_\ell+1}(\xi,y),0).
		\eean
		In particular
		\be\label{e:representation_H_y}
		H(\xi,y) = \kappa \sum^{\infty}_{\ell=1}\gamma^{\sigma^0_\ell}\, S(B_1^{-\sigma^0_\ell+1}(\xi,y),0).
		\ee
	\end{pr}
	
	We now combine the representation formulas for $S$ obtained in Section~\ref{s:SBR} with Proposition~\ref{p:representation_H_explicit} to derive a more explicit representation of $H$. In what follows, we work with the second identity in Proposition~\ref{p:representation_H_explicit}, which is well suited to our subsequent analysis. To prepare for the study of local times of \emph{tilted} (or \emph{drifted}) Takagi curves, we derive a representation of
	$$
	H(\xi,y)+\beta yS(\xi,0),
	$$
	where $\beta\in\mathbb{R}$ is arbitrary.

	
	\begin{co}\label{c:representation_simple_H}
		Let $\xi, y\in[0,1].$ Let $X_n = 1_{\{\by_n=0\}}-1_{\{\by_n=1\}}, R_n = 1_{\{\by_n=1\}}$ for $n\ge 1$. Then we have
		\be\label{e:representation_simple_H_y}
		H(\xi,y) + \beta y S(\xi,0)=
		\frac{\kappa^2}{2} \Big[\sum_{n=1}^\infty \gamma^n X_n - \sum_{n=1}^\infty \sum_{p=1}^\infty \gamma^n (\eh)^{p}\, X_{p+n}\,X_n\Big] + (1+\beta) y S(\xi,0).
		\ee
	\end{co}
	
	\pf According to Proposition \ref{p:representation_H_explicit} 
	we have to calculate $S(B_1^{-\sigma^0_\ell+1}(\xi,y),0)$ for $\xi, y \in[0,1], \ell\in\N$ using the formulas provided by Corollaries \ref{c:simple_representation_S}, \ref{c:simple_representation_S_oneterm}. To exploit them, we first note that the dyadic sequence associated with $B_1^{-\sigma^0_\ell+1}(\xi, y)$ is given by\\ $(\cdots,\bxi_{-1},\bxi_0,\by_1,\by_2,\cdots, \by_{\sigma^0_{\ell}-2},\by_{\sigma^0_\ell-1}).$  We therefore can write
	\bea
	S(B_1^{-\sigma^0_\ell+1}(\xi, y),0) &=& \sum_{m=1}^{\sigma^0_\ell-1} \kappa^{m+1} \big[1_{\{\by_{\sigma^0_\ell-m}=0\}}- 1_{\{\by_{\sigma^0_\ell-m}=1\}}\big]+ \kappa^{\sigma^0_\ell-1}\, S(\xi,0)\\
	&=& \kappa^{\sigma^0_\ell}\,\sum_{n=1}^{\sigma^0_\ell-1} \kappa^{-n+1}\,X_n + \kappa^{\sigma^0_\ell-1} S(\xi,0).
	\eea
	
	Now insert this formula into the representation of Proposition \ref{p:representation_H}, and observe that $\gamma \kappa = \eh$ and that $R_\ell = \frac{1-X_\ell}{2}$ to obtain
	
	\bea
	&&\hspace{.5cm}H(\xi,y) + \beta y S(\xi,0)\\
	&&\hspace{.3cm}= \sum_{\ell=1}^\infty \gamma^{\sigma^0_\ell}\Big[ \kappa^{\sigma^0_\ell+1}\,\sum_{n=1}^{\sigma^0_\ell-1} \kappa^{-n+1}\,X_n + \kappa^{\sigma^0_\ell} S(\xi,0)\Big] + \beta y S(\xi,0)\\
	&&\hspace{.3cm} = \frac{\kappa^2}{2} \sum_{\ell=1}^\infty (\eh)^{\ell}\, (1-X_\ell)\, \sum_{n=1}^{\ell-1} \kappa^{-n} X_n + (1+\beta) y S(\xi,0) \\
	&&\hspace{.3cm} = \frac{\kappa^2}{2} \sum_{n=1}^\infty \sum_{\ell=n+1}^\infty (\eh)^{\ell}\, (1-X_\ell)\,\kappa^{-n} X_n + (1+\beta) y S(\xi,0) \\
	&&\hspace{.3cm} = \frac{\kappa^2}{2} \sum_{n=1}^\infty \sum_{p=1}^\infty (\eh)^{p+n}\, (1-X_{p+n})\,\kappa^{-n} X_n + (1+\beta) y S(\xi,0) \\
	&&\hspace{.3cm} = \frac{\kappa^2}{2} \Big[\sum_{n=1}^\infty \gamma^n X_n - \sum_{n=1}^\infty \sum_{p=1}^\infty \gamma^n (\eh)^{p}\, X_{p+n}\,X_n\Big] + (1+\beta) y S(\xi,0). 
	\eea
	
	This gives the desired representation. \epf

	\subsection{Moments of iterated integrals in Rademacher calculus}\label{ss:moments}

	The representation of $H$ obtained in Corollary~\ref{c:representation_simple_H} reveals an unexpected connection with Malliavin calculus on non-Gaussian probability spaces. Indeed, up to elementary deterministic terms, the quantity
	\[
	H(\xi,\cdot)+\beta\cdot S(\xi,0)
	\]
	is composed of constant multiples of
	\[
	\sum_{n\ge1}\gamma^nX_n
	\quad\text{and}\quad
	\sum_{n\ge1}\sum_{m\ge1}
	1_{{n<m}}
	\kappa^{-n}2^{-m}X_nX_m,
	\]
	where $(X_n)_{n\ge1}$ is a sequence of independent symmetric Bernoulli random variables, also known as a \emph{Rademacher sequence}. These are precisely multiple stochastic integrals belonging to the first and second Rademacher chaoses. The corresponding Malliavin calculus, which shares many features with the classical Gaussian theory, was developed by Nourdin, Peccati and Reinert \cite{nourdinpeccatireinert10} in connection with probabilistic and statistical applications.
	
	Our approach to the existence of occupation densities requires estimates of exponential moments involving linear combinations of first- and second-order Rademacher chaoses. To this end, we introduce the kernels
	\[
	f(i)=\gamma^i,
	\qquad
	h(i,j)=\kappa^{-i}2^{-j}1_{{i<j}},
	\qquad i,j\in\mathbb N,
	\]
	together with the associated multiple stochastic integrals
	\[
	I_1(f)=\sum_{n\ge1}f(n)X_n,
	\qquad
	I_2(h)=\sum_{n\ge1}\sum_{m\ge1}h(n,m)X_nX_m.
	\]
	Since multiple stochastic integrals are most naturally formulated in terms of symmetric kernels, we begin by deriving a representation of the symmetrization of $h$.

	\begin{lem}\label{l:symmetrization:g}
		Let $g$ be the symmetrization of $h$, i.e. $g(i,j) = \eh (h(i,j) + h(j,i)), i,j\in \N.$ Then, as in classical Malliavin calculus, $I_2(g) = I_2(h)$, and for $i,j\in\N$ we have
		$$ g(i,j) = \eh\, \gamma^{\frac{i+j}{2}}\, (\frac{\kappa}{2})^{\frac{|i-j|}{2}}\, 1_{\{i\not= j\}}.$$
	\end{lem}
	
	\pf
	We may write
	\bea
	g(i,j) &=& \eh (\kappa^{-i} 2^{-j}\, 1_{\{i<j\}} + \kappa^{-j} 2^{-i} 1_{\{j<i\}})\\
	&=& \eh \kappa^{-i\wedge j} 2^{-i\vee j} 1_{\{i\not= j\}} = \eh \kappa^{-i\wedge j} (\gamma \kappa)^{i\vee j} 1_{\{i\not= j\}}\\
	&=& \eh \gamma^{i\vee j}\, \kappa^{i\vee j-i\wedge j}\, 1_{\{i\not= j\}} = \eh \gamma^{i\vee j}\, \kappa^{|i-j|}\, 1_{\{i\not= j\}},
	\eea
	
	and alternatively
	\bea
	g(i,j) &=& \eh  \kappa^{-i\wedge j} 2^{-i\vee j} 1_{\{i\not= j\}} = \eh (2\gamma)^{i\wedge j} 2^{-i\vee j} 1_{\{i\not= j\}}\\
	&=& \eh \gamma^{i\wedge j} 2^{-(i\vee j - i\wedge j)} 1_{\{i\not= j\}} = \eh \gamma^{i\wedge j} 2^{-|i-j|} 1_{\{i\not= j\}}.
	\eea
	Taking roots of the alternative expressions obtained, we finally have
	\bea
	g(i,j) &=& \eh \gamma^{\frac{i\vee j}{2}} \kappa^{\frac{|i-j|}{2}}\,\,\gamma^{\frac{i\wedge j}{2}} 2^{-\frac{|i-j|}{2}} 1_{\{i\not= j\}}\\
	&=& \eh \gamma^{\frac{i+j}{2}} (\frac{\kappa}{2})^{\frac{|i-j|}{2}} 1_{\{i\not= j\}}.
	\eea
	This is the desired formula. \epf
	
	In order to calculate the Fourier transform of (constants multiplied to $f$ resp. $g$ are unessential)
	$$\R\ni u\mapsto \exp(i u (I_1(f)+ I_2(g))),$$
	we will mainly be occupied with computing the moments
	$$E(I_1(f)^n I_2(g)^m)$$
	for $n,m\ge 0.$ In fact, it will be sufficient to confine our attention to even $n$, since the transformation $X \mapsto -X$ for the sequence $X = (X_i)_{i\in\N}$ is measure preserving, and via this transformation $I_1(f)$ maps to $- I_1(f)$, while $I_2(g)$ maps to $I_2(g)$. So for odd $n$, we get $E(I_1(f)^n I_2(g)^m) = - E(I_1(f)^n I_2(g)^m),$ proving that
	$$E(I_1(f)^n I_2(g)^m) = 0\quad\mbox{for odd}\quad n\in\N.$$
	
	For the calculus of the moments, let us briefly recall some notation well known from combinatorics, and encountered for instance in Malliavins calculus. We denote the set of all nonnegative integer valued sequences for which almost all components vanish by $S_0$. For $p = (p_i)_{i\in\N} \in S_0$ we write
	$$|p| = \sum_{i\ge 1} p_i,\quad p! = \prod_{i\ge 1} p_i!, \quad \binom{|p|}{p} = \frac{|p|!}{p!} = \frac{|p|!}{\prod_{i\ge 1} p_i!}.$$
	
	\subsubsection{Moments of  Rademacher integrals in the first chaos}\par\medskip
	
	Let us first suppose $m=0$ and compute $E(I_1(f)^{2n})$ for $n\ge 0.$ We may simply write
	$$E(I_1(f)^{2n}) = \sum_{j_1,\cdots, j_{2n}=1}^\infty \gamma^{\sum_{i=1}^{2n} j_i} E(\prod_{i=1}^{2n} X_{j_i}).$$
	Now for $k\in\N$ and a $2n$ tuple of indices $(j_1,\cdot, j_{2n})$ let $T_k = \{i: 1\le i\le 2n, j_i=k\}.$ Then $(T_k)_{k\in \N}$ is a partition of $\{1,\cdots,2n\}$, and for almost all $k\in\N$ we have $T_k=\emptyset.$ So we may write
	$$ E(\prod_{i=1}^{2n} X_{j_i}) = E(\prod_{k\ge 1} X_k^{|T_k|}),$$
	and from the definition of our Bernoulli variables it is immediately clear that
	$ E(\prod_{k\ge 1} X_k^{|T_k|}) =1,$
	if $|T_k|$ is even for any $k\in\N$, while $ E(\prod_{k\ge 1} X_k^{|T_k|}) =0$ if this is not the case. Consequently we have
	\bean\label{e:firstchaos-development}
	&&\hspace{2cm} E(I_1^{2n}(f))\\
	\nonumber &&\hspace{1.5cm} = \sum_{(T_i)_{i\ge 1}\,\mbox{partition of}\,\{1,\cdots,2n\}, |T_i|\,\mbox{even}}\prod_{k\ge 1} \gamma^{k |T_k|}\\
	\nonumber &&\hspace{1.5cm} =\sum_{p\in S_0\, \mbox{even},\, |p|=2n} |\{(T_i)_{i\in\N}\,\mbox{partition of}\,\{1,\cdots,2n\}, |T_i|= p_i\,\mbox{for}\, i\ge 1\}| \prod_{k\ge 1} \gamma^{k p_k}\\
	\nonumber &&\hspace{1.5cm} =\sum_{p\in S_0\, \mbox{even},\, |p|=2n} \binom{2n}{p} \prod_{k\ge 1} \gamma^{k p_k}.
	\eean
	For the last equation in \eqref{e:firstchaos-development} we use the well known fact from combinatorics that the number of partitions of a set of $2n$ elements into subsets of $p = (p_i)_{i\ge 1}$ elements is given by the multinomial coefficient $\binom{2n}{p}$ for multi-indices $p\in S_0$ such that $|p|=2n$. We can summarize our findings in the following statement.
	
	\begin{pr}\label{p:momentsfirstchaos}
		For $n\ge 1$ we have
		$$E(I_1^{2n}(f)) = \sum_{p\in S_0\, \mbox{even},\, |p|=2n} \binom{2n}{p} \prod_{k\ge 1} \gamma^{k p_k}.$$
	\end{pr}
	
	Let us finally use the result of Proposition \ref{p:momentsfirstchaos} to countercheck the \emph{dual} of the well known result from Section \ref{s:SBR}. We obtain
	\bean\label{e:countercheck}
	E(\cos(u I_1(f))) &=& \sum_{n=0}^\infty \frac{(-1)^n u^{2n}}{(2n)!} E(I_1(f)^{2n})\\
	\nonumber &=& \sum_{n=0}^\infty \frac{(-1)^n u^{2n}}{(2n)!} \sum_{p\in S_0\, \mbox{even},\, |p|=2n} \binom{2n}{p} \prod_{k\ge 1} \gamma^{k p_k}\\
	\nonumber&=& \sum_{p\in S_0\, \mbox{even}}\frac{(-1)^{\frac{|p|}{2}}}{|p|!}\binom{|p|}{p} \prod_{k\ge 1} (u \gamma^k)^{p_k}\\
	\nonumber &=& \sum_{p\in S_0\, \mbox{even}}\frac{(-1)^{\frac{|p|}{2}}}{p!} \prod_{k\ge 1} (u \gamma^k)^{p_k}\\
	\nonumber &=& \sum_{p\in S_0\, \mbox{even}} \prod_{k\ge 1}\frac{(-1)^{\frac{p_k}{2}}}{p_k!} (u \gamma^k)^{p_k}= \prod_{k\in\N} \cos(u \gamma^k).
	\eean
	
	This formula could of course have been obtained more directly, without the detour of the calculation of even moments of $I_1(f)$ by taking into account the independence of the $X_i, i\ge 1,$ and elementary trigonometric identities.
	
	\subsubsection{Moments of  Rademacher integrals in the second chaos}\label{sss:moments_secondchaos}
	
	Let us next assume $n=0$ and compute $E(I_2(h)^m)$ for $m\ge 0$ by a similar line of reasoning. We first write
	\begin{align}\label{eqI2m}
		E(I_2(g)^{m}) = (\frac{1}{2})^m\sum_{j_1,\cdots, j_{2m}=1}^\infty 1_{\{j_1\not=j_2,\cdots, j_{2m-1}\not= j_{2m}\}} \gamma^{\frac{\sum_{i=1}^{2m} j_i}{2}} \prod_{r=1}^m (\frac{\kappa}{2})^{\frac{|j_{2r-1}-j_{2r}|}{2}} E(\prod_{i=1}^{2m} X_{j_i}).
	\end{align}
	We therefore now have to reduce our space of admissible $2m$-tuples of integers $(j_1,\cdots,j_{2m})$ to the set
	$$\Delta_m = \{(j_1,j_2,\cdots, j_{2m-1}, j_{2m}): j_{2i-1}\not= j_{2i}\,\mbox{for all}\, 1\le i\le m\}.$$
	For $k\in\N$ and a $2m$-tuple of indices $(j_1,\cdots,j_{2m})\in \Delta_m$ let
	$$T_k = \{i: 1\le i\le 2m, j_i=k\}.$$
	Again, $(T_k)_{k\in \N}$ is a partition of $\{1,\cdots,2m\}$, and for almost all $k\in\N$ we have $T_k=\emptyset.$ Moreover, all partition sets belong to
	$$D_m = \{ T: T \subset \{1,\cdots, 2m\}, \{2i-1, 2i\}\not\subset T\,\mbox{for all}\, 1\le i\le m\}.$$
	For a partition $(T_k)_{k\in\N}$ in $D_m$ we have, as before,
	$$ E(\prod_{i=1}^{2m} X_{j_i}) = E(\prod_{k\ge 1} X_k^{|T_k|}),$$
	and from the definition of our Bernoulli variables it is immediately clear that
	$ E(\prod_{k\ge 1} X_k^{|T_k|}) =1,$
	if $|T_k|$ is even for any $k\in\N$, while $ E(\prod_{k\ge 1} X_k^{|T_k|}) =0$ if this is not the case. To describe the moments, we need one more notation. For $h, k\in \N, h\not= k,$ let
	$$q_{hk} = |\{l:1\le l\le m, \{2l-1,2l\}\subset T_h\cup T_k\}|.$$
	In these terms, we may write
	\bean\label{eqI2m1} 
	&&\hspace{.4cm} E(I_2(g)^m)\\
	&&\hspace{.2cm} = (\eh)^m \sum_{(T_i)_{i\ge 1}\subset D_m\,\mbox{even partition of}\,\{1,\cdots,2m\}}\prod_{k\ge 1} \gamma^{\frac{k}{2} |T_k|}\cdot \prod_{h < k} (\frac{\kappa}{2})^{\frac{k-h}{2} q_{hk}} \notag \\
	&&\hspace{.2cm} = (\eh)^m \sum_{p\in S_0\, \mbox{even},\, |p|=2m}\prod_{k\ge 1} \gamma^{\frac{k}{2} p_k} \notag \\
	&&\hspace{3cm}\cdot \sum_{(T_i)_{i\ge 1}\subset D_m\,\mbox{even partition of}\,\{1,\cdots,2m, |T_i|=p_i\,\mbox{for}\,i\in\N\}} \prod_{h < k} (\frac{\kappa}{2})^{\frac{k-h}{2} q_{hk}}.\notag 
	\eean
	We denote the interaction term appearing in the last line of the preceding inequality by
	$$L_m(p) = \sum_{(T_i)_{i\ge 1}\subset D_m\,\mbox{even partition of}\,\{1,\cdots,2m\}, |T_i|=p_i\,\mbox{for}\,i\in\N} \prod_{h < k} (\frac{\kappa}{2})^{\frac{k-h}{2} q_{hk}},$$
	for an even multi-index $p\in S_0$ such that $|p|=2m.$

	To further compute the function $L_m(p)$, even for $p\in S_0$ not necessarily even, it will be simpler to argue for the complement of $\Delta_m$ in the set if all tuples $(j_1,\cdots, j_{2m})$ of integers. For this purpose, for $H\subset \{1,\cdots,m\}$ let
	$$\Delta_m(H) = \{(j_i,\cdots,j_{2m}): j_{2l-1}=j_{2l}\,\mbox{for}\,l\in H,\quad j_{2l-1}\not= j_{2l}\,\mbox{for}\, l\notin H\},$$
	and correspondingly
	$$D_m(H) = \{ T: T \subset \{1,\cdots, 2m\}, \{2l-1, 2l\}\subset T\,\mbox{for}\,l\in H,\,\{2l-1, 2l\}\not\subset T\,\mbox{for}\,l\not\in H\}.$$
	
	Then, by definition, $(\Delta_m(H))_{H\subset \{1,\cdots,m\}}$ is a pairwise disjoint partition of the set of all $(j_1,\cdots, j_{2m})\in \N^{2m}$, and we may identify $D_m(\emptyset) = D_m$. Therefore we can write
	\bean\label{e:estimate_cardinality_1}
	&&\hspace{0.5cm} L_m(p)\\
	\nonumber &&\hspace{.3cm} = \sum_{(T_i)_{i\ge 1}\subset D_m\,\mbox{partition of}\,\{1,\cdots,2m\}, |T_i|=p_i\,\mbox{for}\,i\in\N} \prod_{h < k} (\frac{\kappa}{2})^{\frac{k-h}{2} q_{hk}}\\
	\nonumber&&\hspace{.3cm} = \sum_{(T_i)_{i\ge 1}\,\mbox{partition of}\,\{1,\cdots,2m\}, |T_i|=p_i\,\mbox{for}\,i\in\N} \prod_{h < k} (\frac{\kappa}{2})^{\frac{k-h}{2} q_{hk}}\\
	\nonumber &&\hspace{.6cm} - \sum_{\emptyset\not= H\subset \{1,\cdots, 2m\}} \sum_{(T_i)_{i\ge 1}\subset D_m(H)\,\mbox{partition of}\,\{1,\cdots,2m\}, |T_i|=p_i\,\mbox{for}\,i\in\N} \prod_{h < k} (\frac{\kappa}{2})^{\frac{k-h}{2} q_{hk}}.
	\eean
	Now define
	$$N_m(p) = \sum_{(T_i)_{i\ge 1}\,\mbox{partition of}\,\{1,\cdots,2m\}, |T_i|=p_i\,\mbox{for}\,i\in\N} \prod_{h < k} (\frac{\kappa}{2})^{\frac{k-h}{2} q_{hk}}.$$
	Now note that the quantities appearing in the last line of \eqref{e:estimate_cardinality_1} are invariant for permutation of index pairs $(j_{2l-1},j_{2l})$. Hence we may abbreviate for $0\le s\le m$ $D_m(s) = D_m(\{1,\cdots,2s\})$ and have in continuation of \eqref{e:estimate_cardinality_1}
	
	\bean\label{e:estimate_cardinality_2}
	&& L_m(p) = N_m(p)\\
	\nonumber&& \hspace{1.5cm}- \sum_{1\le s\le m} \binom{m}{s} \sum_{(T_i)_{i\ge 1}\subset D_m(s)\,\mbox{partition of}\,\{1,\cdots,2m\}, |T_i|=p_i\,\mbox{for}\,i\in\N} \prod_{h < k} (\frac{\kappa}{2})^{\frac{k-h}{2} q_{hk}}.
	\eean
	
	Recall that partitions $(T_l)_{l\ge 1}$ for some $(j_1,\cdots, j_{2m})\in \N_0^{\N}$ just describe the sets
	$$T_l = \{i: j_i=l\}.$$ In this setting, for a partition $(T_l)_{l\ge 1}$ in $D_m(s)$ let us define $\pi^1_s((j_1,\cdots, j_{2m})) = (j_1,\cdots, j_{2s})$, and $\pi_s^2((j_1,\cdots, j_{2m})) = (j_{2s+1},\cdots, j_{2m}),$
	$$U_l(s) = \{i: j_i=l,\,\,\, i\in\pi^1(D_m(s))\},\quad V_l(s) = \{i: j_i=l,\,\,\, i\in\pi^2(D_m(s))\}.$$
	In these terms, we may identify further, with $\tilde{q}_{hk}$ defined as $q_{hk}$, but referring to the partition $(V_l(s))_{l\ge 1}$
	
	\bean\label{e:estimate_cardinality_4}
	&&\hspace{3.5cm} L_m(p)\\
	\nonumber&& \hspace{.2cm}= N_m(p)- \sum_{1\le s\le m} \binom{m}{s}\,\sum_{q+r=p, q\,\mbox{even}} \sum_{\underset{(U_l(s))_{l\ge 1}\,\mbox{partition in}\, \pi^1(D_m(s)),|U_l(s)|=q_l,} {(V_l(s))_{l\ge 1}\,\mbox{partition in}\, \pi^2(D_m(s)), |V_l(s)|=r_l}}\prod_{h < k} (\frac{\kappa}{2})^{\frac{k-h}{2} \tilde{q}_{hk}}\\
	\nonumber&&\hspace{.1cm} = N_m(p) - \sum_{1\le s\le m} \binom{m}{s}\,\sum_{q+r=p, |q|=2s,\, q\,\mbox{even}}
	\frac{1}{2^s} \frac{s!}{(\frac{q}{2})!} \sum_{(V_l(s)_{l\ge 1}\,\mbox{partition in}\, \pi^2(D_m(s)), |V_l(s)|=r_l}\prod_{h < k} (\frac{\kappa}{2})^{\frac{k-h}{2} \tilde{q}_{hk}}\\
	\nonumber&&\hspace{.1cm} = N_m(p) - \sum_{1\le s\le m} \binom{m}{s}\,\sum_{q+r=p, |q|=2s,\, q\,\mbox{even}}
	\frac{1}{2^s} \frac{s!}{(\frac{q}{2})!} L_{m-s}(r)\\
	&&\nonumber\hspace{.2cm} = N_m(p) - \sum_{1\le \frac{|q|}{2}\le m,\, q\,\mbox{even}}
	\frac{1}{2^{\frac{|q|}{2}}} \frac{m!}{(m-\frac{|q|}{2})!\,(\frac{q}{2})!} L_{m-\frac{|q|}{2}}(p-q).
	\eean
	
	From \eqref{e:estimate_cardinality_4} it becomes clear that we now can argue recursively to compute $L_m(p)$. Indeed, writing $q_1, r_1$ instead of $q, r$ we obtain in the next step
	\bean\label{e:estimate_cardinality_5}
	&&\hspace{3cm} N_m(p) - \sum_{1\le \frac{|q_1|}{2}\le m,\, q_1\,\mbox{even}}
	\frac{1}{2^{\frac{|q_1|}{2}}} \frac{m!}{(m-\frac{|q_1|}{2})!\,(\frac{q_1}{2})!} L_{m-\frac{|q_1|}{2}}(p-q_1)\\
	\nonumber&&\hspace{.2cm} = N_m(p) - \sum_{1\le \frac{|q_1|}{2}\le m,\, q_1\,\mbox{even}}
	\frac{1}{2^{\frac{|q_1|}{2}}} \frac{m!}{(m-\frac{|q_1|}{2})!\,(\frac{q_1}{2})!}\times  \big[ N_{m-\frac{|q_1|}{2}}(p-q_1) \\
	&&\nonumber \hspace{1cm}- \sum_{1\le \frac{|q_2|}{2}\le m-\frac{|q_1|}{2}}
	\frac{1}{2^{\frac{|q_2|}{2}}} \frac{(m-\frac{|q_1|}{2})!}{(m-\frac{|q_1|}{2}-\frac{|q_2|}{2})!\,(\frac{q_2}{2})!} L_{m-\frac{|q_1|}{2}-\frac{|q_2|}{2}}(p-q_1-q_2)\big]\\
	\nonumber&&\hspace{.2cm} = N_m(p) - \sum_{1\le \frac{|q_1|}{2}\le m,\, q_1\,\mbox{even}}
	\frac{1}{2^{\frac{|q_1|}{2}}} \frac{m!}{(m-\frac{|q_1|}{2})!\,(\frac{q_1}{2})!} N_{m-\frac{|q_1|}{2}}(p-q_1)\\
	\nonumber&&\hspace{1cm} + \sum_{2\le \frac{|q_1|}{2}+\frac{|q_2|}{2}\le m,\, q_1, q_2\,\mbox{even}}
	\frac{1}{2^{\frac{|q_1|}{2}+\frac{|q_2|}{2}}} \frac{m!}{(m-\frac{|q_1|}{2}-\frac{|q_2|}{2})!\,(\frac{q_1}{2})!\,(\frac{q_2}{2})!} L_{m-\frac{|q_1|}{2}-\frac{|q_2|}{2}}(p-q_1-q_2).
	\eean
	
	Since we interpret $\frac{|q_1|}{2}$ as cardinality of partitions of $\{1,\cdots,\frac{|q_1|}{2}\}$, we may accordingly consider $\frac{|q_2|}{2}$ as cardinality of partitions of $\{\frac{|q_1|}{2}+1,\cdots,\frac{|q_1|}{2}+\frac{|q_2|}{2}\}$. For convenience, we take the union of the two families of disjoint partitions, and simply write $\frac{|q_2|}{2}$ for the cardinality of the unions. This means that we replace $\frac{|q_1|}{2}+\frac{|q_2|}{2}$ simply by $\frac{|q_2|}{2}$, and the multi-index $q_1+q_2$ by $q_2$. With these modifications, \eqref{e:estimate_cardinality_5} becomes
	
	\bean\label{e:estimate_cardinality_6}
	&&\hspace{3cm} N_m(p) - \sum_{1\le \frac{|q_1|}{2}\le m,\, q_1\,\mbox{even}}
	\frac{1}{2^{\frac{|q_1|}{2}}} \frac{m!}{(m-\frac{|q_1|}{2})!\,(\frac{q_1}{2})!} N_{m-\frac{|q_1|}{2}}(p-q_1)\\
	\nonumber&&\hspace{1cm} + \sum_{2\le \frac{|q_1|}{2}+\frac{|q_2|}{2}\le m,\, q_1, q_2\,\mbox{even}}
	\frac{1}{2^{\frac{|q_1|}{2}+\frac{|q_2|}{2}}} \frac{m!}{(m-\frac{|q_1|}{2}-\frac{|q_2|}{2})!\,(\frac{q_1}{2})!\,(\frac{q_2}{2})!} L_{m-\frac{|q_1|}{2}-\frac{|q_2|}{2}}(p-q_1-q_2)\\
	\nonumber&&=N_m(p) - \sum_{1\le \frac{|q_1|}{2}\le m,\, q_1\,\mbox{even}}
	\frac{1}{2^{\frac{|q_1|}{2}}} \frac{m!}{(m-\frac{|q_1|}{2})!\,(\frac{q_1}{2})!} N_{m-\frac{|q_1|}{2}}(p-q_1)\\
	\nonumber&&\hspace{1cm} + \sum_{2\le \frac{|q_2|}{2}\le m,\, q_2\,\mbox{even}}
	\frac{1}{2^{\frac{|q_2|}{2}}} \frac{m!}{(m-\frac{|q_2|}{2})!\,(\frac{q_2}{2})!} L_{m-\frac{|q_2|}{2}}(p-q_2)
	\eean
	
	After at most $m+1$ iterations we get an expression composed of terms, the $l$th one of which for $l\ge 1$ reads
	\bean\label{e:estimate_cardinality_7}
	&&(-1)^l \sum_{l\le \frac{|q_l|}{2}\le m,\, q_1\,\mbox{even}}
	\frac{1}{2^{\frac{|q_l|}{2}}} \frac{m!}{(m-\frac{|q_l|}{2})!\,(\frac{q_l}{2})!} N_{m-\frac{|q_l|}{2}}(p-q_l)
	\eean
	For $l=0$ the expression in \eqref{e:estimate_cardinality_5} evidently reduces to $N_m(p)$.
	Therefore we obtain
	
	\begin{pr}\label{p:cardinality_L}
		For $p\in \N_0^{\N}$ such that $|p|=2m$ let
		$$L_m(p)
		= \sum_{(T_i)_{i\ge 1}\subset D_m\,\mbox{partition of}\,\{1,\cdots,2m\}, |T_i|=p_i\,\mbox{for}\,i\in\N} \prod_{h < k} (\frac{\kappa}{2})^{\frac{k-h}{2} q_{hk}},$$
		and
		$$N_m(p) = \sum_{(T_i)_{i\ge 1}\,\mbox{partition of}\,\{1,\cdots,2m\}, |T_i|=p_i\,\mbox{for}\,i\in\N} \prod_{h < k} (\frac{\kappa}{2})^{\frac{k-h}{2} q_{hk}}.$$
		Then we have
		\bean\label{estimate_cardinality_6}
		&&\hspace{0.4cm}L_m(p) \\
		\nonumber&&=  \sum_{l=0}^m (-1)^l \sum_{l\le \frac{|q_l|}{2}\le m,\, q_1\,\mbox{even}}
		\frac{1}{2^{\frac{|q_l|}{2}}} \frac{m!}{(m-\frac{|q_l|}{2})!\,(\frac{q_l}{2})!} N_{m-\frac{|q_l|}{2}}(p-q_l)\\
		\nonumber&& = \sum_{0\le \frac{|q|}{2}\le m ,\, q\,\mbox{even},\,|q|\in 4\N} \frac{1}{2^{\frac{|q|}{2}}} \frac{m!}{(m-\frac{|q|}{2})!\,(\frac{q}{2})!} N_{m-\frac{|q|}{2}}(p-q).
		\eean
	\end{pr}
	
	\pf
	The last line of \eqref{e:estimate_cardinality_6} follows from the cancellation of terms for which $\frac{|q|}{2}$ is odd, as seen from the previous line. \epf
	
	It remains to compute $N_m(p)$ for partitions $(T_i)_{i\ge 1}$ of $\{1,\cdots,2m\}$, $m\in\N.$ This will turn out much simpler than the direct computation of the corresponding $L_m(p).$ In fact, for such a partition we let for $i\in\N$
	$$ U_i = T_i \cap \{2,4,\cdots, 2m\},\quad V_i = T_i\cap \{1,3,\cdots, 2m-1\}.$$
	Then $|U_i|+|V_i|=|T_i| = p_i, i\in\N$, and $|U_i|, |V_i|\le m.$ In these terms, by definition, we have
	$$q_{kh} = |(U_k-1) \cap V_h| + |(U_h-1) \cap V_k|,\quad k,h\in\N.$$
	But this may be simplified further. By writing $(j_1,j_3,\cdots,j_{2m-1}), (j_2,j_4,\cdots, j_{2m})$\\ as $(f_1,\cdots,f_m), (g_1,\cdots,g_m)$ and sampling $(U_i)_{i\ge 1}$ from the tuples $(f_1,\cdots, f_m)$ and $(V_i)_{i\ge 1}$ from $(g_1,\cdots, g_m)$, re-interpreting
	$q_{kh} = |U_k\cap V_h|+|U_h\cap V_k|$, we can rewrite
	\bean\label{estimate_cardinality_6}
	&&\hspace{1.5cm}N_m(p)\\
	\nonumber&&\hspace{1cm} = \sum_{u+v=p} \sum_{(U_i)_{i\ge 1},(V_i)_{i\ge 1} \,\mbox{partitions of}\,\{1,\cdots,m\}, |U_i|=u_i, |V_i|=v_i} \prod_{h < k} (\frac{\kappa}{2})^{\frac{k-h}{2} q_{hk}}.
	\eean
	Now let $r_{kh} = |U_k\cap V_h|.$ Then $q_{kh} = r_{kh}+r_{hk}, h, k\in\N.$ We know that $\sum_{h\in\N} r_{kh}= u_k, \sum_{k\in\N} r_{kh}=v_h.$ The cardinality of the set of partitions $(U_i)_{i\ge 1}$ of $\{1,\cdots, m\}$ is given by $\binom{m}{u}$. Now each one of the sets of this partition $U_k$ gets subdivided into sets of cardinalities $r_{kh}, h\in\N$. The cardinalities of these sets of subdivisions are given by $\binom{u_k}{(r_{kh})_{h\in\N}}$. Hence we obtain a total of
	$$\sum_{\sum_{h\in\N} (r_{kh}+r_{hk})=p_k, k\in\N}\binom{m}{u}\cdot\prod_{k\in\N}\binom{u_k}{(r_{kh})_{h\in\N}} = \sum_{\sum_{h\in\N} (r_{kh}+r_{hk})=p_k, k\in\N} \binom{m}{(r_{kh})_{k,h\in\N}}$$
	configurations. So we finally obtain
	$$N_m(p) = \sum_{\sum_{h\in\N} q_{kh}=p_k, k\in\N}  \binom{m}{(r_{kh})_{k,h\in\N}} \, \prod_{h < k} (\frac{\kappa}{2})^{\frac{k-h}{2} q_{kh}}.$$
	
	We may summarize our findings in
	
	\begin{pr}\label{p:cardinality_N}
		For $p\in \N_0^{\N}$ such that $|p|=2m$ let
		$$N_m(p) = \sum_{(T_i)_{i\ge 1}\,\mbox{partition of}\,\{1,\cdots,2m\}, |T_i|=p_i\,\mbox{for}\,i\in\N} \prod_{h < k} (\frac{\kappa}{2})^{\frac{k-h}{2} q_{hk}}.$$
		Then we have
		\be\label{estimate_cardinality_7}
		N_m(p)=
		\sum_{\sum_{h\in\N} q_{kh}=p_k, k\in\N}  \binom{m}{(r_{kh})_{k,h\in\N}} \, \prod_{h < k} (\frac{\kappa}{2})^{\frac{k-h}{2} q_{kh}}.
		\ee
	\end{pr}

	The following summarizes our findings from Propositions \ref{p:cardinality_L} and \ref{p:cardinality_N} about the moments of $I_2(g)$.
	
	\begin{pr}\label{p:momentssecondschaos}
		For $m\ge 1$ we have
		\bean\label{e:momentssecondchaos}
		&&\hspace{3cm} E(I_2(g)^m) = \\
		\nonumber &&\sum_{p,q,r\in S_0\, p,q\mbox{even},\, |p|=2m, |q|\in 4\N, \sum_{h\in\N} (r_{\cdot h}+r_{h\cdot})+q=p}
		\prod_{k\ge 1} \gamma^{\frac{k}{2} p_k}\cdot (\eh)^{\frac{|p|}{2}+\frac{|q|}{2}} \binom{m}{\frac{q}{2}\,\, (r_{kh})_{k,h\in\N}} \cdot \prod_{h < k} (\frac{\kappa}{2})^{\frac{k-h}{2} (r_{kh}+r_{hk})}.
		\eean
	\end{pr}
	
	\subsubsection{Moments of products of Rademacher integrals in the first and second chaos}
	
	Let us finally compute the moments of $I_1(f)^{2n} I_2(g)^m$ for $n,m\in\N$. According to the cases discussed above we are allowed to write
	\bea
	&&\hspace{.5cm} I_1(f)^{2n} I_2(g)^m\\
	&&\hspace{.2cm} = (\eh)^m \sum_{j_1,\cdots, j_{2n}, j_{2n+1,\cdots,j_{2n+2m}=1}}^\infty \gamma^{\sum_{i=1}^{2n} j_i+\frac{\sum_{l=1}^{2m} j_{2n+l}}{2}} 1_{\{j_{2n+1}\not=j_{2n+2},\cdots, j_{2n+2m-1}\not= j_{2n+2m}\}}\\
	&&\hspace{2cm}\cdot\prod_{r=1}^m (\frac{\kappa}{2})^\frac{|j_{2n+2r-1}-j_{2n+2r}|}{2} E(\prod_{i=1}^{2n+2m} X_{j_i}).
	\eea
	In natural extensions of notation used above let
	$$\Delta_{n,m} = \{(j_1,\cdots,j_{2n},j_{2n+1},\cdots, j_{2n+2m}): j_{2n+2i-1}\not= j_{2n+2i}\,\mbox{for all}\, 1\le i\le m\},$$
	for $(j_1,\cdots, j_{2n+2m})\in \Delta_{n,m}$ and $k\in\N$ let
	$$T_k = \{i: 1\le i\le 2n+2m, j_i=k\},$$
	$$D_{n,m} = \{ T: T \subset \{1,\cdots, 2n+2m\}, \{2n+2i-1, 2n+2i\}\not\subset T\,\mbox{for all}\, 1\le i\le m\},$$
	and finally for $h,k\in\N, h\not= k$ let
	$$q_{kh} = |\{l:1\le l\le m, \{2n+2l-1,2n+2l\}\subset T_h\cup T_k\}|.$$
	
	Again, $(T_k)_{k\ge 1}$ is a partition of $\{1,\cdots, 2n+2m\}$, and the condition that $E(\prod_{i=1}^{2n+2m} X_{j_i})$ either vanish or be 1 leads to the conclusion that only \emph{even} partitions $(T_k)_{k\ge 1}$ of $\{1,\cdots,2n+2m\}$ give nonzero contributions. So we obtain the formula
	\bean\label{e:firstandsecondchaos-development}
	&&\hspace{.4cm} E(I_1(f)^{2n} I_2(g)^m)\\
	\nonumber &&\hspace{.2cm} = (\eh)^m \sum_{(T_i)_{i\ge 1}\subset D_{n,m}\,\mbox{even partition of}\,\{1,\cdots,2n+2m\}}\\
	\nonumber &&\hspace{3cm}\cdot \prod_{k\ge 1} \gamma^{k |T_k\cap\{1,\cdots,2n\}|+\frac{k}{2} |T_k\cap\{2n+1,\cdots,2n+2m\}|}
	\cdot \prod_{h < k} (\frac{\kappa}{2})^{\frac{k-h}{2} q_{hk}}.
	\eean
	
	\eqref{e:firstandsecondchaos-development} stipulates the use of separate multi-indices to denote partitions of $\{1,\cdots, 2n\}$ and $\{2n+1,\cdots, 2n+2m\}.$ We note
	$$U_i = T_i\cap \{1,\cdots,2n\},\quad V_i = T_i\cap \{2n+1,\cdots,2n+2m\},\quad i\in\N.$$ Then we can give a multi-index version of \eqref{e:firstandsecondchaos-development} by
	\bean\label{e:firstandsecondchaos-multiindex}
	&&\hspace{.5cm} E(I_1(f)^{2n} I_2(g)^m)\\
	\nonumber &&\hspace{.3cm} = \sum_{q,r\in S_0, q+r\,\mbox{even}, |q|=2n, |r|=2m} (\eh)^{\frac{|r|}{2}}\cdot \prod_{k\ge 1} \gamma^{k q_k+\frac{k}{2}r_k}\\
	\nonumber &&\hspace{1cm}  \sum_{(U_i)_{i\in\N}\,\mbox{part. of}\, \{1,\cdots,2n\}, (V_i-2n)_{i\in\N}\subset D_m\,\mbox{part. of}\, \{1,\cdots,2m\}, |U_i|=q_i, |V_i|=r_i} \prod_{h < k} (\frac{\kappa}{2})^{\frac{k-h}{2} q_{kh}}\\
	\nonumber &&\hspace{.3cm} = \sum_{q,r\in S_0, q+r\,\mbox{even}, |q|=2n, |r|=2m} (\eh)^{\frac{|r|}{2}}\cdot \prod_{k\ge 1} \gamma^{k q_k+\frac{k}{2}r_k}\\
	\nonumber &&\hspace{1cm} \binom{2n}{q}\,\sum_{(V_i-2n)_{i\in\N}\subset D_m\,\mbox{part. of}\, \{1,\cdots,2m\}, |V_i|=r_i} \prod_{h < k} (\frac{\kappa}{2})^{\frac{k-h}{2} q_{kh}}\\
	\nonumber &&\hspace{.3cm}= \sum_{q,r\in S_0, q+r\,\mbox{even}, |q|=2n, |r|=2m} (\eh)^{\frac{|r|}{2}}\cdot \prod_{k\ge 1} \gamma^{k q_k+\frac{k}{2}r_k} \binom{2n}{q}\, L_m(r).
	\eean
	
	The simplification leading to the second equation follows from the combinatorial arguments given in the discussion of the moments of $I_1(f)^{2n}$. To identify\\ $\sum_{(V_i-2n)_{i\in\N}\subset D_m\,\mbox{part. of}\, \{1,\cdots,2m\}, |V_i|=r_i} \prod_{h < k} (\frac{\kappa}{2})^{\frac{k-h}{2} q_{kh}}$ with $L(r)$ (without the condition that the related partitions be even) already defined in the discussion of the moments of $I_2(g)^m$, we just have to realize that taking intersections with $\{2n+1,\cdots,2n+2m\}$ and shifting by $2n$ produces a partition sequence of $\{1,\cdots,2m\}$ in $D_m$, if we start with a partition sequence of $\{1,\cdots, 2n+2m\} \subset D_{n,m}$. Moreover, $q_{kh}$ defined above correspond exactly to the coefficients defined in the discussion of the moments of $I_2(g)^m$, if we take those with respect to the partition sequence of $\{2n+1,\cdots,2n+2m\}$ shifted by $2n$. In summary we obtain
	
	\begin{pr}\label{p:momentsfirstandsecondschaos}
		For $n,m\in \N$ we have
		$$E(I_1(f)^{2n} I_2(g)^m) = \sum_{q,r\in S_0, q+r\,\mbox{even}, |q|=2n, |r|=2m} (\eh)^{\frac{|r|}{2}}\cdot \prod_{k\ge 1} \gamma^{k q_k+\frac{k}{2}r_k} \binom{2n}{q}\, L_m(r).$$
	\end{pr}
	
	Note that while $q+r$ in the formula of Lemma \ref{p:momentsfirstandsecondschaos} is an even partition, $q$ and $r$ need not be even. This can be seen by the following Example.\par
	
	{\bf Ex. 1}: Let $m=3, n=2$, and take the following 10-tuple of indices:
	$$j_1=j_5=j_7=j_9=1, j_2=j_6=2, j_3=j_8=3, j_4=j_{10}=4.$$
	This tuple obviously belongs to $D_{2,3}$.
	Then we have $$T_1=\{1,5,7,9\}, T_2 = \{2,6\}, T_3 = \{3,8\}, T_4 = \{4, 10\}, T_k=\emptyset\, \mbox{for}\, k\ge 5.$$
	And so $(T_k)_{k\ge 1}$ is even, while
	$$U_1 = T_1\cap \{1,\cdots,4\} = \{1\}, V_1 = T_1\cap \{5,\cdots, 10\} = \{5,7,9\},$$
	hence neither $(U_k)_{k\ge 1}$ nor $(V_k)_{k\ge 1}$ is even. \par
	
	Using the results of subsection \ref{sss:moments_secondchaos}, and recalling that the $L_m(r)$ of Proposition \ref{p:momentsfirstandsecondschaos} is a simple extension of its analogue in subsection \ref{sss:moments_secondchaos} we can give a still more explicit expression for the formula in Proposition \ref{p:momentsfirstandsecondschaos}.
	
	\begin{pr}\label{p:momentsfirstandsecondchaos_explicit}
		For $n,m\in \N$ we have
		\bean\label{e:momentsfirstandsecondchaos_explicit}
		&&\hspace{2cm} E(I_1(f)^{2n} I_2(g)^m) =\\
		\nonumber && \sum_{o,p,q,r\in S_0, o+p, q\,\mbox{even}, |o|=2n, |p|=2m, |q|\in 4\N, \sum_{h\in\N} (r_{\cdot h}+r_{h\cdot})+q=p}\\
		\nonumber&&\hspace{2cm} (\eh)^{\frac{|p|}{2}+\frac{|q|}{2}}\cdot \prod_{k\ge 1} \gamma^{k o_k+\frac{k}{2}p_k} \binom{2n}{o}\,  \binom{m}{\frac{q}{2}\,\, (r_{kh})_{k,h\in\N}} \cdot \prod_{h < k} (\frac{\kappa}{2})^{\frac{k-h}{2} (r_{kh}+r_{hk})}.
		\eean
	\end{pr}
	
	\pf This is a combination of Propositions \ref{p:momentssecondschaos} and \ref{p:momentsfirstandsecondschaos}. \epf

	\subsection{The Fourier transform}\label{ss:fouriertransform}
	
	Let us finally apply Proposition \ref{p:momentsfirstandsecondchaos_explicit} to compute the Fourier transforms of our occupation time functional, recalling from Corollary \ref{c:representation_simple_H} that we have to deal with linear combinations of components from the first and second Rademacher chaos. We obtain
	
	\begin{pr}\label{p:FT-with-rademacher}
		For $u\in\R$, and $a, b\in\R$ we have
		\bea
		&&\hspace{.7cm} E[\exp\big\{iu( aI_1(f)+ bI_2(g))\big\}]\\
		&&\hspace{.5cm}= \sum_{o,q,r\in S_0, o+\sum_{h\in\N} (r_{\cdot h}+r_{h\cdot}), q\,\mbox{even}, |q|\in 4\N}\\
		\nonumber&&\hspace{2cm} (-1)^{\frac{|o|}{2}} \prod_{k\ge 1} \frac{(u^2a^2\gamma^{2k})^{\frac{o_k}{2}}}{o_k!} \,\,(-1)^{\frac{|q|}{4}}\prod_{k\ge 1}  \frac{(\frac{u^2b^2}{16}\gamma^{2k})^{\frac{q_k}{4}}}{(\frac{q_k}{2})!}\,(-1)^{\frac{|r|}{2}} \prod_{h< k}\frac{(\frac{u^2b^2}{16}\gamma^{2k} (\frac{\kappa}{2})^{2(k-h)})^{\frac{r_{hk}+r_{kh}}{4}}}{(r_{kk})! (r_{kh})!(r_{hk})!}.
		\eea
		
	\end{pr}
	
	\pf
	Indeed, for $u\in\R$ the Proposition \ref{p:momentsfirstandsecondchaos_explicit} 
	yields
	\bean\label{e:FT-with-rademacher}
	&&\hspace{2cm} E[\exp\big\{iu( aI_1(f)+ bI_2(g))\big\}]\\
	\nonumber &&\hspace{.2cm} = \sum_{n,m=0}^\infty \frac{(iu)^{n+m}a^nb^m}{(n)! (m)!}  E(I_1(f)^{n} I_2(g)^{m})\\
	\nonumber &&\hspace{.2cm} = \sum_{n,m=0}^\infty \frac{(-1)^{\frac{n+m}{2}}}{(n)! (m)!} u^{n+m}a^nb^m\cdot \sum_{o,p,q,r\in S_0, o+p, q\,\mbox{even}, |o|=n, |p|=2m, |q|\in 4\N, \sum_{h\in\N} (r_{\cdot h}+r_{h\cdot})+q=p}\\
	\nonumber&&\hspace{2cm} (\eh)^{\frac{|p|}{2}+\frac{|q|}{2}}\cdot \prod_{k\ge 1} \gamma^{k o_k+\frac{k}{2}p_k} \binom{n}{o}\,  \binom{m}{\frac{q}{2}\,\, (r_{kh})_{k,h\in\N}} \cdot \prod_{h < k} (\frac{\kappa}{2})^{\frac{k-h}{2} (r_{kh}+r_{hk})}\\
	\nonumber &&\hspace{.2cm}= \sum_{o,p,q,r\in S_0, o+p, q\,\mbox{even}, |q|\in 4\N, \sum_{h\in\N} (r_{\cdot h}+r_{h\cdot})+q=p}\\
	\nonumber&&\hspace{2cm} (\eh)^{\frac{|p|}{2}+\frac{|q|}{2}}\cdot (-1)^{\frac{|o|}{2}} \prod_{k\ge 1} \frac{(u^2a^2\gamma^{2k})^{\frac{o_k}{2}}}{o_k!} \,\,(-1)^{\frac{|p|}{4}}\prod_{k\ge 1}  \frac{(u^2b^2\gamma^{2k})^{\frac{p_k}{4}}}{(\frac{q_k}{2})! (r_{k\cdot})!}   \cdot \prod_{h < k} (\frac{\kappa}{2})^{\frac{k-h}{2} (r_{kh}+r_{hk})}.
	\eean
	Now we use the equality $\frac{|p|}{2}= \frac{|q|}{2}+|r|$, to distribute the appearing powers of $\eh$ and $-1$ to the last two products. We remark that given $q$ is even, the condition $o+p$ even readily translates into the condition that $o+\sum_{h\in\N} (r_{\cdot h}+r_{h\cdot})$ is even.
	Finally note that the seemingly awkward notation in the products $\prod_{k\ge 1}\cdots$ is necessary, since $\frac{o_k}{2}$ and $\frac{p_k}{4}$ need not be integer. This is the desired formula. \epf
	
	The formula produced by Proposition \ref{p:FT-with-rademacher} can be given a much simpler and concise form. To derive it, we shall start giving a simpler formula for the central part in the following Lemma.
	
	\begin{lem}\label{l:FT-centralpart}
		We have for $b\in\R$
		$$\sum_{q\,\mbox{even},\, |q|\in 4\N} (-1)^{\frac{|q|}{4}}\prod_{k\ge 1}  \frac{(\frac{u^2b^2}{16}\gamma^{2k})^{\frac{q_k}{4}}}{(\frac{q_k}{2})!} = \cos(\frac{ub}{4} \frac{\gamma}{1-\gamma}).$$
	\end{lem}
	
	\pf Using $i^2 = -1$, we may rewrite
	\bea
	&&\sum_{q\,\mbox{even},\, |q|\in 4\N} (-1)^{\frac{|q|}{4}}\prod_{k\ge 1}  \frac{(\frac{u^2b^2}{16}\gamma^{2k})^{\frac{q_k}{4}}}{(\frac{q_k}{2})!}
	= \sum_{q\,\mbox{even},\, |q|\in 4\N} \prod_{k\ge 1}  \frac{(i \frac{ub}{4}\gamma^{k})^{\frac{q_k}{2}}}{(\frac{q_k}{2})!}\\
	&&\hspace{.5cm} = \sum_{p\in S_0,\, |p|\in 2\N} \prod_{k\ge 1}  \frac{(i \frac{ub}{4}\gamma^{k})^{p_k}}{p_k!}
	= \lim_{N\to\infty} \sum_{p_1,\cdots, p_N\ge 0, |p|\in 2\N} \prod_{k = 1}^N  \frac{(i \frac{ub}{4}\gamma^{k})^{p_k}}{p_k!}.
	\eea
	
	Now fix $N\in\N$. We can write
	\bea
	&&\sum_{p_1,\cdots, p_N\ge 0, |p|\in 2\N} \prod_{k = 1}^N  \frac{(i \frac{ub}{4}\gamma^{k})^{p_k}}{p_k!}\\
	&&\hspace{.5cm} = \sum_{p_1,\cdots, p_{N-1}\ge 0, |p|\in 2\N} \prod_{k = 1}^{N-1}  \frac{(i \frac{ub}{4}\gamma^{k})^{p_k}}{p_k!}\\
	&&\hspace{3cm} \cdot \big[\sum_{k=0}^\infty \frac{(i \frac{ub}{4}\gamma^{N})^{2k}}{(2k)!} 1_{2\N}(p_1+\cdots+p_{N-1}) + \sum_{k=0}^\infty \frac{(i \frac{ub}{4}\gamma^{N})^{2k+1}}{(2k+1)!} 1_{2\N+1}(p_1+\cdots+p_{N-1})\big]\\
	&&\hspace{.5cm} = \sum_{p_1,\cdots, p_{N-1}\ge 0, |p|\in 2\N} \prod_{k = 1}^{N-1}  \frac{(i \frac{ub}{4}\gamma^{k})^{p_k}}{p_k!}\\
	&&\hspace{3cm} \cdot \big[\cos(\frac{ub}{4} \gamma^N) 1_{2\N}(p_1+\cdots+p_{N-1}) + i \sin(\frac{ub}{4} \gamma^N) 1_{2\N+1}(p_1+\cdots+p_{N-1})\big]\\
	&&\hspace{.5cm} = \sum_{p_1,\cdots, p_{N-1}\ge 0, |p|\in 2\N} \prod_{k = 1}^{N-1}  \frac{(i \frac{ub}{4}\gamma^{k})^{p_k}}{p_k!}\\
	&&\hspace{3cm} \cdot \eh \big[\exp(i \frac{ub}{4} \gamma^N) + (-1)^{p_1+\cdots+p_{N-1}} \exp(-i \frac{ub}{4} \gamma^N) \big]\\
	&&\hspace{.5cm} = \eh \big[\sum_{p_1,\cdots, p_{N-1}\ge 0, |p|\in 2\N} \prod_{k = 1}^{N-1}  \frac{(i \frac{ub}{4}\gamma^{k})^{p_k}}{p_k!}\,\exp(i \frac{ub}{4} \gamma^N)\\
	&&\hspace{3cm} + \sum_{p_1,\cdots, p_{N-1}\ge 0, |p|\in 2\N} \prod_{k = 1}^{N-1}  \frac{(-i \frac{ub}{4}\gamma^{k})^{p_k}}{p_k!} \exp(-i \frac{ub}{4} \gamma^N) \big]\\
	&&\hspace{.5cm} = \eh\big[ \exp(i \frac{ub}{4} \sum_{k=1}^N \gamma^k) + \exp(-i \frac{ub}{4} \sum_{k=1}^N \gamma^k)\big]\\
	&&\hspace{.5cm} = \cos(\frac{ub}{4} \sum_{k=1}^N \gamma^k).
	\eea
	
	It remains to let $N\to\infty$ to obtain the desired result. \epf
	
	The evaluation of the remainder of the expression provided by Proposition \ref{p:FT-with-rademacher} is done in the following Lemma. Arguments are similar to the ones used in the preceding proposition.
	
	\begin{pr}\label{II:FT-externalparts}
		For $u,a,b\in \mathbb{R}_{+}$ we get%
		\begin{eqnarray}\label{eqfi1}
			\phi (u) &=&E\left[ \exp (iu(aI_{1}(f)+bI_{2}(g)))\right]  \\
			&=&\cos (\frac{ub}{4}\frac{\gamma }{1-\gamma })\times \lim_{N\longrightarrow
				\infty }\mathcal{H}_{N},  \notag
		\end{eqnarray}
		\bigskip where%
		\begin{eqnarray*}
			\mathcal{H}_{N} &=&\frac{1}{2^{N}}\sum_{m=1}^{N}\sum_{1\leq
				k_{1}<...<k_{m}\leq N}\exp \Big\{ -i\sum_{j=1}^{m}ua\gamma
			^{k_{j}}+i\sum_{k\notin \left\{ k_{1},...,k_{m}\right\} }ua\gamma ^{k}
			\\
			&&+i\sum_{1\leq j\leq m,1\leq l\leq m}(\frac{ub}{4}\gamma ^{k_{j}}(\frac{%
				\kappa }{2})^{\left\vert k_{j}-k_{l}\right\vert })^{\frac{1}{2}%
			}-i\dsum\limits_{1\leq j\leq m,h\notin \left\{ k_{1},...,k_{m}\right\} }(%
			\frac{ub}{4}\gamma ^{k_{j}}(\frac{\kappa }{2})^{\left\vert
				k_{j}-h\right\vert })^{\frac{1}{2}} \\
			&& -i\sum_{k\notin \left\{ k_{1},...,k_{m}\right\} }\sum_{j=1}^{m}(%
			\frac{ub}{4}\gamma ^{k}(\frac{\kappa }{2})^{\left\vert k-k_{j}\right\vert
			})^{\frac{1}{2}})+i\sum_{k,h\notin \left\{ k_{1},...,k_{m}\right\} }(\frac{ub%
			}{4}\gamma ^{k}(\frac{\kappa }{2})^{\left\vert k-h\right\vert })^{\frac{1}{2}%
			})\Big\} \text{.}
		\end{eqnarray*}
		
	\end{pr}
	
	\pf
	Using Lemma \ref{p:FT-with-rademacher}, we have that%
	\begin{eqnarray*}
		&&\sum_{o,r\in S_{0},o+\sum_{h\in \mathbb{N}}(r_{\cdot h}+r_{h\cdot })\text{
				even}}(i)^{\left\vert o\right\vert }\dprod\limits_{k\geq 1}\frac{%
			(u^{2}a^{2}\gamma ^{2k})^{\frac{o_{k}}{2}}}{o_{k}!}(i)^{\left\vert
			r\right\vert }\dprod\limits_{h,k\geq 1}\frac{(\frac{ub}{4}\gamma ^{k}(\frac{%
				\kappa }{2})^{\left\vert k-h\right\vert })^{\frac{r_{kh}}{2}}}{(r_{kh})!} \\
		&=&\sum_{o,r\in S_{0},o+\sum_{h\in \mathbb{N}}(r_{\cdot h}+r_{h\cdot })\text{
				even}}\dprod\limits_{k\geq 1}\frac{(iua\gamma ^{k})^{o_{k}}}{o_{k}!}%
		\dprod\limits_{h,k\geq 1}\frac{(i(\frac{ub}{4}\gamma ^{k}(\frac{\kappa }{2}%
			)^{\left\vert k-h\right\vert })^{\frac{1}{2}})^{r_{kh}}}{(r_{kh})!} \\
		&=&\lim_{N\longrightarrow \infty }\sum_{(o_{j})_{1\leq j\leq
				N},(r_{i_{1}j_{1}})_{1\leq i_{1},j_{1}\leq N}:o+\sum_{h\in \mathbb{N}%
			}(r_{\cdot h}+r_{h\cdot })\text{ even}}\dprod\limits_{k=1}^{N}(\frac{%
			(iua\gamma ^{k})^{o_{k}}}{o_{k}!})\dprod\limits_{h,k= 1}^N\frac{(i(\frac{ub%
			}{4}\gamma ^{k}(\frac{\kappa }{2})^{\left\vert k-h\right\vert })^{\frac{1}{2}%
			})^{r_{kh}}}{(r_{kh})!} \\
		&=&\lim_{N\longrightarrow \infty }\sum_{(o_{j})_{1\leq j\leq
				N},(r_{i_{1}j_{1}})_{1\leq i_{1},j_{1}\leq N}}\dprod\limits_{k=1}^{N}(\frac{1%
		}{2}(1+(-1)^{o_{k}+\sum_{h=1}^{N}(r_{kh}+r_{hk})}))\times \\
		&&\dprod\limits_{k=1}^{N}(\frac{(iua\gamma ^{k})^{o_{k}}}{o_{k}!}%
		)\dprod\limits_{h,k= 1}^N\frac{(i(\frac{ub}{4}\gamma ^{k}(\frac{\kappa }{2}%
			)^{\left\vert k-h\right\vert })^{\frac{1}{2}})^{r_{kh}}}{(r_{kh})!} \\
		&=&\lim_{N\longrightarrow \infty }\sum_{(o_{j})_{1\leq j\leq
				N},(r_{i_{1}j_{1}})_{1\leq i_{1},j_{1}\leq N}}(\frac{1}{2}%
		)^{N}(1+\sum_{m=1}^{N}\sum_{1\leq k_{1}<...<k_{m}\leq
			N}\dprod\limits_{j=1}^{m}(-1)^{o_{k_{j}}+%
			\sum_{h=1}^{N}(r_{k_{j}h}+r_{hk_{j}})})\times \\
		&&\dprod\limits_{k=1}^{N}(\frac{(iua\gamma ^{k})^{o_{k}}}{o_{k}!}%
		)\dprod\limits_{h,k= 1}^N\frac{(i(\frac{ub}{4}\gamma ^{k}(\frac{\kappa }{2}%
			)^{\left\vert k-h\right\vert })^{\frac{1}{2}})^{r_{kh}}}{(r_{kh})!} \\
		&=&\lim_{N\longrightarrow \infty }((\frac{1}{2})^{N}\sum_{(o_{j})_{1\leq
				j\leq N},(r_{i_{1}j_{1}})_{1\leq i_{1},j_{1}\leq N}}\dprod\limits_{k=1}^{N}(%
		\frac{(iua\gamma ^{k})^{o_{k}}}{o_{k}!})\dprod\limits_{h,k= 1}^N\frac{(i(%
			\frac{ub}{4}\gamma ^{k}(\frac{\kappa }{2})^{\left\vert k-h\right\vert })^{%
				\frac{1}{2}})^{r_{kh}}}{(r_{kh})!} \\
		&&+(\frac{1}{2})^{N}\sum_{m=1}^{N}\sum_{1\leq k_{1}<...<k_{m}\leq
			N}\sum_{(o_{j})_{1\leq j\leq N},(r_{i_{1}j_{1}})_{1\leq i_{1},j_{1}\leq
				N}}\dprod\limits_{j=1}^{m}(-1)^{o_{k_{j}}+%
			\sum_{h=1}^{N}(r_{k_{j}h}+r_{hk_{j}})}\times \\
		&&\dprod\limits_{k=1}^{N}(\frac{(iua\gamma ^{k})^{o_{k}}}{o_{k}!}%
		)\dprod\limits_{h,k= 1}^N\frac{(i(\frac{ub}{4}\gamma ^{k}(\frac{\kappa }{2}%
			)^{\left\vert k-h\right\vert })^{\frac{1}{2}})^{r_{kh}}}{(r_{kh})!})\text{.}
	\end{eqnarray*}
	
	Further, we see that%
	\begin{eqnarray*}
		&&\sum_{(o_{j})_{1\leq j\leq N},(r_{i_{1}j_{1}})_{1\leq i_{1},j_{1}\leq
				N}}\dprod\limits_{j=1}^{m}(-1)^{o_{k_{j}}+%
			\sum_{h=1}^{N}(r_{k_{j}h}+r_{hk_{j}})}\dprod\limits_{k=1}^{N}(\frac{%
			(iua\gamma ^{k})^{o_{k}}}{o_{k}!})\dprod\limits_{h,k= 1}^N\frac{(i(\frac{ub%
			}{4}\gamma ^{k}(\frac{\kappa }{2})^{\left\vert k-h\right\vert })^{\frac{1}{2}%
			})^{r_{kh}}}{(r_{kh})!} \\
		&=&\sum_{(o_{j})_{1\leq j\leq N},(r_{i_{1}j_{1}})_{1\leq i_{1},j_{1}\leq
				N}}\dprod\limits_{j=1}^{m}(-1)^{\sum_{h=1}^{N}(r_{k_{j}h}+r_{hk_{j}})}\dprod%
		\limits_{j=1}^{m}(\frac{(-iua\gamma ^{k_{j}})^{o_{k_{j}}}}{o_{k_{j}}!}%
		)\dprod\limits_{k\notin \left\{ k_{1},...,k_{m}\right\} }(\frac{(iua\gamma
			^{k})^{o_{k}}}{o_{k}!})\times \\
		&&\dprod\limits_{1\leq j\leq m,1\leq h\leq N}\frac{(i(\frac{ub}{4}\gamma
			^{k_{j}}(\frac{\kappa }{2})^{\left\vert k_{j}-h\right\vert })^{\frac{1}{2}%
			})^{r_{k_{j}h}}}{(r_{k_{j}h})!}\dprod\limits_{(k,h):k\notin \left\{
			k_{1},...,k_{m}\right\} }\frac{(i(\frac{ub}{4}\gamma ^{k}(\frac{\kappa }{2}%
			)^{\left\vert k-h\right\vert })^{\frac{1}{2}})^{r_{kh}}}{(r_{kh})!} \\
		&=&\sum_{(o_{j})_{1\leq j\leq N},(r_{i_{1}j_{1}})_{1\leq i_{1},j_{1}\leq
				N}}\dprod\limits_{j=1}^{m}(-1)^{\sum_{h=1}^{N}r_{hk_{j}}}\dprod%
		\limits_{j=1}^{m}(\frac{(-iua\gamma ^{k_{j}})^{o_{k_{j}}}}{o_{k_{j}}!}%
		)\dprod\limits_{k\notin \left\{ k_{1},...,k_{m}\right\} }(\frac{(iua\gamma
			^{k})^{o_{k}}}{o_{k}!})\times \\
		&&\dprod\limits_{1\leq j\leq m,1\leq h\leq N}\frac{(-i(\frac{ub}{4}\gamma
			^{k_{j}}(\frac{\kappa }{2})^{\left\vert k_{j}-h\right\vert })^{\frac{1}{2}%
			})^{r_{k_{j}h}}}{(r_{k_{j}h})!}\dprod\limits_{(k,h):k\notin \left\{
			k_{1},...,k_{m}\right\} }\frac{(i(\frac{ub}{4}\gamma ^{k}(\frac{\kappa }{2}%
			)^{\left\vert k-h\right\vert })^{\frac{1}{2}})^{r_{kh}}}{(r_{kh})!}\text{.}
	\end{eqnarray*}%
	We have that%
	\begin{eqnarray*}
		\dprod\limits_{j=1}^{m}(-1)^{\sum_{h=1}^{N}r_{hk_{j}}} 	&=&\dprod\limits_{j=1}^{m}\Big\{
		\dprod\limits_{l=1}^{m}(-1)^{r_{k_{l}k_{j}}}\dprod\limits_{k\notin \left\{
			k_{1},...,k_{m}\right\} }(-1)^{r_{kk_{j}}}\Big\} \\
		&=&\dprod\limits_{l=1}^{m}\dprod\limits_{j=1}^{m}(-1)^{r_{k_{j}k_{l}}}\dprod%
		\limits_{j=1}^{m}\dprod\limits_{k\notin \left\{ k_{1},...,k_{m}\right\}
		}(-1)^{r_{kk_{j}}} \\
		&=&\dprod\limits_{j=1}^{m}\dprod\limits_{l=1}^{m}(-1)^{r_{k_{j}k_{l}}}\dprod%
		\limits_{j=1}^{m}\dprod\limits_{k\notin \left\{ k_{1},...,k_{m}\right\}
		}(-1)^{r_{kk_{j}}} \\
		&=&\dprod\limits_{j=1}^{m}\dprod\limits_{l=1}^{m}(-1)^{r_{k_{j}k_{l}}}\dprod%
		\limits_{k\notin \left\{ k_{1},...,k_{m}\right\}
		}\dprod\limits_{j=1}^{m}(-1)^{r_{kk_{j}}}\text{.}
	\end{eqnarray*}%
	Hence%
	\begin{eqnarray*}
		&&\sum_{(o_{j})_{1\leq j\leq N},(r_{i_{1}j_{1}})_{1\leq i_{1},j_{1}\leq
				N}}\dprod\limits_{j=1}^{m}(-1)^{\sum_{h=1}^{N}r_{hk_{j}}}\dprod%
		\limits_{j=1}^{m}(\frac{(-iua\gamma ^{k_{j}})^{o_{k_{j}}}}{o_{k_{j}}!}%
		)\dprod\limits_{k\notin \left\{ k_{1},...,k_{m}\right\} }(\frac{(iua\gamma
			^{k})^{o_{k}}}{o_{k}!})\times \\
		&&\dprod\limits_{1\leq j\leq m,1\leq h\leq N}\frac{(-i(\frac{ub}{4}\gamma
			^{k_{j}}(\frac{\kappa }{2})^{\left\vert k_{j}-h\right\vert })^{\frac{1}{2}%
			})^{r_{k_{j}h}}}{(r_{k_{j}h})!}\dprod\limits_{(k,h):k\notin \left\{
			k_{1},...,k_{m}\right\} }\frac{(i(\frac{ub}{4}\gamma ^{k}(\frac{\kappa }{2}%
			)^{\left\vert k-h\right\vert })^{\frac{1}{2}})^{r_{kh}}}{(r_{kh})!} \\
		&=&\sum_{(o_{j})_{1\leq j\leq N},(r_{i_{1}j_{1}})_{1\leq i_{1},j_{1}\leq
				N}}\dprod\limits_{j=1}^{m}(\frac{(-iua\gamma ^{k_{j}})^{o_{k_{j}}}}{%
			o_{k_{j}}!})\dprod\limits_{k\notin \left\{ k_{1},...,k_{m}\right\} }(\frac{%
			(iua\gamma ^{k})^{o_{k}}}{o_{k}!})\times \\
		&&\dprod\limits_{1\leq j\leq m,1\leq l\leq m}\frac{(i(\frac{ub}{4}\gamma
			^{k_{j}}(\frac{\kappa }{2})^{\left\vert k_{j}-k_{l}\right\vert })^{\frac{1}{2%
			}})^{r_{k_{j}k_{l}}}}{(r_{k_{j}k_{l}})!}\dprod\limits_{1\leq j\leq m,h\notin
			\left\{ k_{1},...,k_{m}\right\} }\frac{(-i(\frac{ub}{4}\gamma ^{k_{j}}(\frac{%
				\kappa }{2})^{\left\vert k_{j}-h\right\vert })^{\frac{1}{2}})^{r_{k_{j}h}}}{%
			(r_{k_{j}h})!}\times \\
		&&\dprod\limits_{k\notin \left\{ k_{1},...,k_{m}\right\}
		}\dprod\limits_{j=1}^{m}\frac{(-i(\frac{ub}{4}\gamma ^{k}(\frac{\kappa }{2}%
			)^{\left\vert k-k_{j}\right\vert })^{\frac{1}{2}})^{r_{kk_{j}}}}{%
			(r_{kk_{j}})!}\dprod\limits_{k,h\notin \left\{ k_{1},...,k_{m}\right\} }%
		\frac{(i(\frac{ub}{4}\gamma ^{k}(\frac{\kappa }{2})^{\left\vert
				k-h\right\vert })^{\frac{1}{2}})^{r_{kh}}}{(r_{kh})!}\text{.}
	\end{eqnarray*}

	Therefore we obtain that%
	\begin{align*}
		&\sum_{(o_{j})_{1\leq j\leq N},(r_{i_{1}j_{1}})_{1\leq i_{1},j_{1}\leq
				N}}\dprod\limits_{j=1}^{m}(-1)^{\sum_{h=1}^{N}r_{hk_{j}}}\dprod%
		\limits_{j=1}^{m}(\frac{(-iua\gamma ^{k_{j}})^{o_{k_{j}}}}{o_{k_{j}}!}%
		)\dprod\limits_{k\notin \left\{ k_{1},\ldots,k_{m}\right\} }(\frac{(iua\gamma
			^{k})^{o_{k}}}{o_{k}!})\times \\
		&\dprod\limits_{1\leq j\leq m,1\leq h\leq N}\frac{(-i(\frac{ub}{4}\gamma
			^{k_{j}}(\frac{\kappa }{2})^{\left\vert k_{j}-h\right\vert })^{\frac{1}{2}%
			})^{r_{k_{j}h}}}{(r_{k_{j}h})!}\dprod\limits_{(k,h):k\notin \left\{
			k_{1},\ldots,k_{m}\right\} }\frac{(i(\frac{ub}{4}\gamma ^{k}(\frac{\kappa }{2}%
			)^{\left\vert k-h\right\vert })^{\frac{1}{2}})^{r_{kh}}}{(r_{kh})!} \\
		=&\dprod\limits_{j=1}^{m}\exp (-iua\gamma ^{k_{j}})\dprod\limits_{k\notin
			\left\{ k_{1},\ldots,k_{m}\right\} }\exp (iua\gamma ^{k})\times \\
		&\dprod\limits_{1\leq j\leq m,1\leq l\leq m}\exp (i(\frac{ub}{4}\gamma
		^{k_{j}}(\frac{\kappa }{2})^{\left\vert k_{j}-k_{l}\right\vert })^{\frac{1}{2%
		}})\dprod\limits_{1\leq j\leq m,h\notin \left\{ k_{1},\ldots,k_{m}\right\}
		}\exp (-i(\frac{ub}{4}\gamma ^{k_{j}}(\frac{\kappa }{2})^{\left\vert
			k_{j}-h\right\vert })^{\frac{1}{2}})\times \\
		&\dprod\limits_{k\notin \left\{ k_{1},\ldots,k_{m}\right\}
		}\dprod\limits_{j=1}^{m}\exp (-i(\frac{ub}{4}\gamma ^{k}(\frac{\kappa }{2}%
		)^{\left\vert k-k_{j}\right\vert })^{\frac{1}{2}})\dprod\limits_{k,h\notin
			\left\{ k_{1},\ldots,k_{m}\right\} }\exp (i(\frac{ub}{4}\gamma ^{k}(\frac{%
			\kappa }{2})^{\left\vert k-h\right\vert })^{\frac{1}{2}})) \\
		=&\exp \Big\{ -i\sum_{j=1}^{m}ua\gamma ^{k_{j}}+i\sum_{k\notin \left\{
			k_{1},\ldots,k_{m}\right\} }ua\gamma ^{k} \\
		&+i\sum_{1\leq j\leq m,1\leq l\leq m}(\frac{ub}{4}\gamma ^{k_{j}}(\frac{%
			\kappa }{2})^{\left\vert k_{j}-k_{l}\right\vert })^{\frac{1}{2}%
		}-i\dsum\limits_{1\leq j\leq m,h\notin \left\{ k_{1},\ldots,k_{m}\right\} }(%
		\frac{ub}{4}\gamma ^{k_{j}}(\frac{\kappa }{2})^{\left\vert
			k_{j}-h\right\vert })^{\frac{1}{2}} \\
		& -i\sum_{k\notin \left\{ k_{1},\ldots,k_{m}\right\} }\sum_{j=1}^{m}(%
		\frac{ub}{4}\gamma ^{k}(\frac{\kappa }{2})^{\left\vert k-k_{j}\right\vert
		})^{\frac{1}{2}}+i\sum_{k,h\notin \left\{ k_{1},\ldots,k_{m}\right\} }(\frac{ub%
		}{4}\gamma ^{k}(\frac{\kappa }{2})^{\left\vert k-h\right\vert })^{\frac{1}{2}%
		}\Big\} \text{.}
	\end{align*}
	So%
	\begin{align*}
		&\sum_{o,r\in S_{0},o+\sum_{h\in \mathbb{N}}(r_{\cdot h}+r_{h\cdot })\text{
				even}}(i)^{\left\vert o\right\vert }\dprod\limits_{k\geq 1}\frac{%
			(u^{2}a^{2}\gamma ^{2k})^{\frac{o_{k}}{2}}}{o_{k}!}(i)^{\left\vert
			r\right\vert }\dprod\limits_{h,k\geq 1}\frac{(\frac{ub}{4}\gamma ^{k}(\frac{%
				\kappa }{2})^{\left\vert k-h\right\vert })^{\frac{r_{kh}}{2}}}{(r_{kh})!} \\
		=&\lim_{N\longrightarrow \infty }((\frac{1}{2})^{N}\sum_{(o_{j})_{1\leq
				j\leq N},(r_{i_{1}j_{1}})_{1\leq i_{1},j_{1}\leq N}}\dprod\limits_{k=1}^{N}(%
		\frac{(iua\gamma ^{k})^{o_{k}}}{o_{k}!})\dprod\limits_{h,k= 1}^N\frac{(i(%
			\frac{ub}{4}\gamma ^{k}(\frac{\kappa }{2})^{\left\vert k-h\right\vert })^{%
				\frac{1}{2}})^{r_{kh}}}{(r_{kh})!} \\
		&+(\frac{1}{2})^{N}\sum_{m=1}^{N}\sum_{1\leq k_{1}<...<k_{m}\leq
			N}\sum_{(o_{j})_{1\leq j\leq N},(r_{i_{1}j_{1}})_{1\leq i_{1},j_{1}\leq
				N}}\dprod\limits_{j=1}^{m}(-1)^{o_{k_{j}}+%
			\sum_{h=1}^{N}(r_{k_{j}h}+r_{hk_{j}})}\times \\
		&\dprod\limits_{k=1}^{N}(\frac{(iua\gamma ^{k})^{o_{k}}}{o_{k}!}%
		)\dprod\limits_{h,k\geq 1}\frac{(i(\frac{ub}{4}\gamma ^{k}(\frac{\kappa }{2}%
			)^{\left\vert k-h\right\vert })^{\frac{1}{2}})^{r_{kh}}}{(r_{kh})!}) \\
		=&\lim_{N\longrightarrow \infty }((\frac{1}{2})^{N}\sum_{m=1}^{N}\sum_{1%
			\leq k_{1}<...<k_{m}\leq N}\sum_{(o_{j})_{1\leq j\leq
				N},(r_{i_{1}j_{1}})_{1\leq i_{1},j_{1}\leq
				N}}\dprod\limits_{j=1}^{m}(-1)^{o_{k_{j}}+%
			\sum_{h=1}^{N}(r_{k_{j}h}+r_{hk_{j}})})\times \\
		&\dprod\limits_{k=1}^{N}(\frac{(iua\gamma ^{k})^{o_{k}}}{o_{k}!}%
		)\dprod\limits_{h,k\geq 1}\frac{(i(\frac{ub}{4}\gamma ^{k}(\frac{\kappa }{2}%
			)^{\left\vert k-h\right\vert })^{\frac{1}{2}})^{r_{kh}}}{(r_{kh})!}) 
	\end{align*}
	\begin{align*}
		=&\lim_{N\longrightarrow \infty }((\frac{1}{2})^{N}\sum_{m=1}^{N}\sum_{1%
			\leq k_{1}<\ldots<k_{m}\leq N}\exp \Big\{ -i\sum_{j=1}^{m}ua\gamma
		^{k_{j}}+i\sum_{k\notin \left\{ k_{1},\ldots,k_{m}\right\} }ua\gamma ^{k}
		\\
		&+i\sum_{1\leq j\leq m,1\leq l\leq m}(\frac{ub}{4}\gamma ^{k_{j}}(\frac{%
			\kappa }{2})^{\left\vert k_{j}-k_{l}\right\vert })^{\frac{1}{2}%
		}-i\dsum\limits_{1\leq j\leq m,h\notin \left\{ k_{1},\ldots,k_{m}\right\} }(%
		\frac{ub}{4}\gamma ^{k_{j}}(\frac{\kappa }{2})^{\left\vert
			k_{j}-h\right\vert })^{\frac{1}{2}} \\
		& -i\sum_{k\notin \left\{ k_{1},\ldots,k_{m}\right\} }\sum_{j=1}^{m}(%
		\frac{ub}{4}\gamma ^{k}(\frac{\kappa }{2})^{\left\vert k-k_{j}\right\vert
		})^{\frac{1}{2}}+i\sum_{k,h\notin \left\{ k_{1},\ldots,k_{m}\right\} }(\frac{ub%
		}{4}\gamma ^{k}(\frac{\kappa }{2})^{\left\vert k-h\right\vert })^{\frac{1}{2}%
		}\Big\} \text{.}
	\end{align*}%
	Finally, the latter representation, Proposition \ref{p:FT-with-rademacher} and
	Lemma \ref{l:FT-centralpart} give the proof.

	\epf

	\subsection{The existence of local time}\label{ss:wintnercases}
	
	In this subsection, we investigate the existence of local times for the Takagi curve $\mathcal{T}$. Our approach is based on the Fourier-analytic criterion developed in the preceding sections together with the Bernoulli convolution representation of the occupation measure. We begin by proving that drifted Takagi curves admit square-integrable local times for almost every $\gamma\in\left(\frac12,1\right)$.

	\subsubsection{Main Results}
	
	The following is one of the main results of the paper. 

	
	\begin{thm}\label{t:integrability_FT_wintnercases}
		For a.e.$\gamma\in(1/2,1)$, it holds $\phi \in
		L^{2}(\mathbb{R})$.
		Therefore $(y\longmapsto H(\xi ,y)-\kappa yS(\xi ,0))$
		possesses a square integrable occupation density.
	\end{thm}

	\begin{proof}
		
		In this proof, we estimate the $L^2$-norm of $\phi$ in terms of the $L^2$-norm of
		$E\left[\exp\big(iu a I_1(f)\big)\right]$, and then invoke
		\cite[Theorem 1.1]{solomyak95} to complete the proof.

		Let $a_{1}\in (0,1)$ and let $%
		a_{1}^{\ast }:=1-a_{1}$. For this purpose, let us recall from \eqref{eqI2m} that%
		\begin{eqnarray*}
			&&E\left[ (I_{2}(g))^{m}\right]  \\
			&=&(\frac{1}{2})^{m}\sum_{j_{1},\ldots,j_{2m}=1}^{\infty }1_{\left\{ j_{1}\neq
				j_{2},\ldots,j_{2m-1}=j_{2m}\right\} }\gamma ^{\frac{\sum_{i=1}^{2m}j_{i}}{2}%
			}\dprod\limits_{r=1}^{m}(\frac{\kappa }{2})^{\frac{\left\vert
					j_{2r-1}-j_{2r}\right\vert }{2}}E\Big[ \dprod\limits_{i=1}^{2m}X_{j_{i}}%
			\Big]  \\
			&=&(\frac{1}{2})^{m}\sum_{j_{1},\ldots,j_{2m}=1}^{\infty }1_{\left\{ j_{1}\neq
				j_{2},\ldots,j_{2m-1}=j_{2m}\right\} }\gamma ^{\frac{a_{1}\sum_{i=1}^{2m}j_{i}}{%
					2}}\dprod\limits_{r=1}^{m}\gamma ^{\frac{a_{1}^{\ast }(j_{2r-1}+j_{2r})}{2}}(%
			\frac{\kappa }{2})^{\frac{\left\vert j_{2r-1}-j_{2r}\right\vert }{2}}E\Big[
			\dprod\limits_{i=1}^{2m}X_{j_{i}}\Big] \text{.}
		\end{eqnarray*}%
		Then using \eqref{eqI2m1} we find that%
		\begin{eqnarray*}
			E\left[ (I_{2}(g))^{m}\right] 
			&=&(\frac{1}{2})^{m}\sum_{p\in S_{0}\text{ even, }\left\vert p\right\vert
				=2m}\dprod\limits_{k\geq 1}\gamma ^{\frac{a_{1}}{2}kp_{k}}\times  \\
			&&\sum_{(T_{i})_{i\geq 1}\subset D_{m}\text{ even partition of }\left\{
				1,\ldots,2m,\left\vert T_{i}\right\vert =p_{i},i\in \mathbb{N}\right\} }\prod_{h<k}\gamma
			^{\frac{a_{1}^{\ast }}{2}(k+h)q_{hk}}(\frac{\kappa }{2})^{^{\frac{k-h}{2}%
					q_{hk}}}.
		\end{eqnarray*}%
		Further, \eqref{e:momentssecondchaos} entails%
		\begin{eqnarray*}
			&&E\left[ (I_{2}(g))^{m}\right]  \\
			&=&\sum_{p,q,r\in S_{0},p,q\text{ even},\left\vert p\right\vert
				=2m,\left\vert q\right\vert \in 4\mathbb{N},\sum_{h\in \mathbb{N}}(r_{\cdot
					h}+r_{h\cdot })+q=p}\dprod\limits_{k\geq 1}\gamma ^{\frac{a_{1}}{2}kp_{k}}(%
			\frac{1}{2})^{\frac{\left\vert p\right\vert }{2}+\frac{\left\vert
					q\right\vert }{2}}\left( 
			\begin{array}{c}
				m \\ 
				\frac{q}{2}(r_{kh})_{k,h\in \mathbb{N}}%
			\end{array}%
			\right) \times  \\
			&&\dprod\limits_{h<k}(\gamma ^{\frac{a_{1}^{\ast }}{2}(k+h)(r_{kh}+r_{hk})}(%
			\frac{\kappa }{2})^{^{\frac{k-h}{2}(r_{kh}+r_{hk})}}).
		\end{eqnarray*}%
		In addition, following the proof of Proposition \ref{p:FT-with-rademacher}, we see that%
		\begin{eqnarray*}
			&&E\left[ \exp (iu(aI_{1}(f)+bI_{2}(g)))\right]  \\
			&=&\sum_{o,q,r\in S_{0},o+\sum_{h\in \mathbb{N}}(r_{\cdot h}+r_{h\cdot }),q%
				\text{ even},\left\vert q\right\vert \in 4\mathbb{N}}(-1)^{\frac{\left\vert
					o\right\vert }{2}}\dprod\limits_{k\geq 1}\frac{(u^{2}a^{2}\gamma ^{2k})^{%
					\frac{o_{k}}{2}}}{o_{k}!}\times  \\
			&&(-1)^{\frac{\left\vert q\right\vert }{4}}\dprod\limits_{k\geq 1}\frac{(%
				\frac{u^{2}b^{2}}{16}(\gamma ^{a_{1}})^{2k})^{\frac{q_{k}}{4}}}{(\frac{q_{k}%
				}{2})!}\times  \\
			&&(-1)^{\frac{\left\vert r\right\vert }{2}}\dprod\limits_{h<k}\frac{(\frac{%
					u^{2}b^{2}}{16}\gamma ^{2(k+a_{1}^{\ast }h)}(\frac{\kappa }{2}%
				)^{^{2(k-h)}})^{\frac{r_{kh}+r_{hk}}{4}}}{(r_{kk})!(r_{kh})!(r_{hk})!}.
		\end{eqnarray*}%
		Using Lemma \ref{l:FT-centralpart}, we obtain  
		\begin{equation*}
			\sum_{q\text{ even},q\in 4\mathbb{N}}(-1)^{\frac{\left\vert q\right\vert }{4}%
			}\dprod\limits_{k\geq 1}\frac{(\frac{u^{2}b^{2}}{16}(\gamma ^{a_{1}})^{2k})^{%
					\frac{q_{k}}{4}}}{(\frac{q_{k}}{2})!}=\cos (\frac{ub}{4}\frac{\gamma ^{a_{1}}%
			}{1-\gamma ^{a_{1}}}).
		\end{equation*}%
		Then, proceeding as in the proof of Proposition \ref{II:FT-externalparts}, we obtain the following representation of $\phi (u)$. For $u,a,b\in \mathbb{R}_{+}$, we have%
		\begin{align}
			\phi (u) =&E\left[ \exp (iu(aI_{1}(f)+bI_{2}(g)))\right] 
			\label{representation} \\
			=&\cos (\frac{ub}{4}\frac{\gamma ^{a_{1}}}{1-\gamma ^{a_{1}}})\times
			\lim_{N\longrightarrow \infty }\mathcal{H}_{N},  \notag
		\end{align}%
		where%
		\begin{eqnarray*}
			\mathcal{H}_{N} &=&\frac{1}{2^{N}}\sum_{m=1}^{N}\sum_{1\leq
				k_{1}<\ldots<k_{m}\leq N}\exp \Big\{ -i\sum_{j=1}^{m}ua\gamma
			^{k_{j}}+i\sum_{k\notin \left\{ k_{1},\ldots,k_{m}\right\} }ua\gamma
			^{k} \\
			&&+i\sum_{1\leq j\leq m,1\leq l\leq m}(\frac{ub}{4}\gamma
			^{k_{j}+a_{1}^{\ast }k_{l}}(\frac{\kappa }{2})^{\left\vert
				k_{j}-k_{l}\right\vert })^{\frac{1}{2}}-i\dsum\limits_{1\leq j\leq m,h\notin
				\left\{ k_{1},\ldots,k_{m}\right\} }(\frac{ub}{4}\gamma ^{k_{j}+a_{1}^{\ast }h}(%
			\frac{\kappa }{2})^{\left\vert k_{j}-h\right\vert })^{\frac{1}{2}} \\
			&& -i\sum_{k\notin \left\{ k_{1},\ldots,k_{m}\right\} }\sum_{j=1}^{m}(%
			\frac{ub}{4}\gamma ^{k+a_{1}^{\ast }k_{j}}(\frac{\kappa }{2})^{\left\vert
				k-k_{j}\right\vert })^{\frac{1}{2}}+i\sum_{k,h\notin \left\{
				k_{1},\ldots,k_{m}\right\} }(\frac{ub}{4}\gamma ^{k+a_{1}^{\ast }h}(\frac{%
				\kappa }{2})^{\left\vert k-h\right\vert })^{\frac{1}{2}}\Big\} \text{.}
		\end{eqnarray*}%
		The second factor in (\ref{representation}) can be rewritten (by using one
		correction term) and we get that%
		\begin{align*}
			&\phi (u) \\
			=&\cos (\frac{ub}{4}\frac{\gamma ^{a_{1}}}{1-\gamma ^{a_{1}}})\times\lim_{N\longrightarrow \infty }\frac{1}{2^{N}}\sum_{K_{1}=1}^{2}\ldots%
			\sum_{K_{N}=1}^{2} \times  \\
			&\exp \Big(i(ua\sum_{\nu =1}^{N}(-1)^{K_{\nu }}\gamma ^{\nu
			}+\left( \frac{ub}{4}\right) ^{\frac{1}{2}}\sum_{\nu ,\mu
				=1}^{N}(-1)^{K_{\nu }}(-1)^{K_{\mu }}\gamma ^{\frac{\nu+a_{1}^{\ast }\mu }{2}%
			}\left( \frac{\kappa }{2}\right) ^{\left\vert \nu -\mu \right\vert /2})\Big).
		\end{align*}%
		Furthermore, if $X_{i},i\geq 1$ is an $i.i.d.$-sequence of standard Bernoulli random variables, then an induction argument shows that the second factor in the above expression can be reformulated as follows:
		\begin{align}
			\phi (u) =&E\left[ \exp (iu(aI_{1}(f)+bI_{2}(g)))\right]   \notag \\
			=&\cos (\frac{ub}{4}\frac{\gamma ^{a_{1}}}{1-\gamma ^{a_{1}}})\times  
			\notag \\
			&\quad E\Big[ \exp (i(ua\sum_{\nu \geq 1}X_{\nu }\gamma ^{\nu }+\left( \frac{ub}{%
				4}\right) ^{\frac{1}{2}}\sum_{\nu ,\mu \geq 1}X_{\nu }X_{\mu }\gamma ^{\frac{%
					\nu+a_{1}^{\ast }\mu }{2}}\left( \frac{\kappa }{2}\right) ^{\left\vert \nu
				-\mu \right\vert /2}))\Big]   \notag \\
			=&\cos \Big(\frac{ub}{4}\frac{\gamma ^{a_{1}}}{1-\gamma ^{a_{1}}}\Big)E\Big[ \exp (i(uaI_{1}(f)+\Big( \frac{ub}{%
				4}\Big) ^{\frac{1}{2}}I_{2}(g^{(1)})))\Big] ,
			\label{NewRepr}
		\end{align}%
		where the new kernel $g^{(1)}$ is defined as $g^{(1)}(i,j)=\gamma ^{\frac{i+a_{1}^{\ast }j}{2}}\left( \frac{%
			\kappa }{2}\right) ^{\left\vert i-j\right\vert /2}$.
		
		Note that%
		\begin{align*}
			&E\Big[ \exp (i(uaI_{1}(f)+\Big( \frac{ub}{%
				4}\Big) ^{\frac{1}{2}}I_{2}(g^{(1)})))\Big]\\ 
			=&\exp \Big(i\sum_{j\geq 1}\Big( \frac{ub}{%
				4}\Big) ^{\frac{1}{2}}g^{(1)}(j,j)\Big)
			E\Big[ \exp (i(uaI_{1}(f)+\Big( \frac{ub}{%
				4}\Big) ^{\frac{1}{2}}I_{2}(g^{\ast })))\Big] ,
		\end{align*}%
		where $g^{\ast }(i,j)=g^{(1)}1_{\left\{ i\neq j\right\} }$. Therefore, we
		obtain that%
		\begin{equation*}
			\big|\phi (u)\big|\leq \Big| E\Big[ \exp (i(uaI_{1}(f)+\Big( \frac{ub}{4}%
			\Big) ^{1/2}I_{2}(g^{\ast })))\Big] \Big|
		\end{equation*}%
		for all $u>0$ and $a_{1}\in (0,1)$. We have that%
		\begin{equation*}
			g^{(1)}(i,j)=\gamma ^{\frac{i+a_{1}^{\ast }j}{2}}\left( \frac{\kappa }{2}%
			\right) ^{\left\vert i-j\right\vert /2}\underset{a_{1}\searrow 0}{%
				\longrightarrow }\gamma ^{\frac{i+j}{2}}\left( \frac{\kappa }{2}\right)
			^{\left\vert i-j\right\vert /2}\text{.}
		\end{equation*}%
		Then it follows from dominated convergence and the definition of the kernel $%
		g$ that%
		\begin{eqnarray}
			\big|\phi (u)\big|&\leq &\underline{\lim }_{a_{1}\longrightarrow 0^{+}}\Big| E%
			\Big[ \exp (i(uaI_{1}(f)+\left( \frac{ub}{4}\right) ^{1/2}I_{2}(g^{\ast })))%
			\Big] \Big|  \notag \\
			&=&\Big|E\Big[ \exp (i(uaI_{1}(f)+\left( ub\right) ^{1/2}I_{2}(g)))%
			\Big] \Big|  \label{IneqR}
		\end{eqnarray}%
		for all $u>0$. Define now%
		\begin{equation*}
			u^{(n)}=\left( u\right) ^{1/2^{n}},b^{(n)}=b^{1/2^{n}}\text{ and }%
			a^{(n)}=\left( u\right) ^{-1/2^{n}}ua\text{.}
		\end{equation*}%
		So repeated application of \eqref{IneqR} implies%
		\begin{align*}
			\big|\phi (u)\big|\leq &\Big| E\Big[ \exp (i(uaI_{1}(f)+( ub)
			^{1/2}I_{2}(g)))\Big] \Big| \\
			=&\Big| E\Big[ \exp (iu^{(1)}(a^{(1)}I_{1}(f)+b^{(1)}I_{2}(g)))\Big]
			\Big| \\
			\leq &\Big| E\Big[ \exp
			(i(u^{(1)}a^{(1)}I_{1}(f)+(u^{(1)}b^{(1)})^{1/2}I_{2}(g)))\Big] \Big|
			\\
			=&\Big| E\Big[ \exp (i(uaI_{1}(f)+(u^{(1)}b^{(1)})^{1/2}I_{2}(g)))%
			\Big] \Big| \\
			=&\Big| E\Big[ \exp (iu^{(2)}(a^{(2)}I_{1}(f)+b^{(2)}I_{2}(g)))\Big]
			\Big|\\
			&\vdots 
		\end{align*}
		\begin{align*}
			&\vdots \\
			\leq &\Big| E\Big[ \exp (iu^{(n)}(a^{(n)}I_{1}(f)+b^{(n)}I_{2}(g)))%
			\Big] \Big| \\
			=&\Big| E\Big[ \exp (i(uaI_{1}(f)+u^{(n)}b^{(n)}I_{2}(g)))\Big]
			\Big| .
		\end{align*}%
		Since $u^{(n)}b^{(n)}\underset{n\longrightarrow \infty }{\longrightarrow }1$
		for all $u,b>0$, we find that%
		\begin{equation*}
			\big|\phi (u)\big|\leq \left\vert E\left[ \exp (i(uaI_{1}(f)+I_{2}(g)))\right]
			\right\vert .
		\end{equation*}%
		Let $\Xi =aI_{1}(f)$ and $\Phi =I_{2}(g)$. Hence,%
		\begin{equation*}
			\left\vert \phi (u)\right\vert ^{2}\leq \left\vert E\left[ \exp (i(u\Xi
			+\Phi ))\right] \right\vert ^{2}
		\end{equation*}%
		for all $u>0$. On the other hand, we have that%
		\begin{equation*}
			E\left[ \exp (i(u\Xi +\Phi ))\right] =E\left[ \exp (iu\Xi )g(\Xi )\right] ,
		\end{equation*}%
		where $g(x):=E	\big[ \exp (i\Phi )	\big| \Xi =x	\big] $ with%
		\begin{equation*}
			\left\vert g(x)\right\vert \leq 1\text{ a.e. on }\left\{ x:f_{\Xi
			}(x)>0\right\}
		\end{equation*}%
		for the square integrable probability density $f_{\Xi }$ of $\Xi $ (which
		exists under our assumptions on $\gamma $). So we can write%
		\begin{equation*}
			E\left[ \exp (i(u\Xi +\Phi ))\right] =\int_{\mathbb{R}}\exp (iux)g(x)f_{\Xi
			}(x)\mathrm{d}x.
		\end{equation*}%
		Thus, we obtain from Plancherel`s theorem that%
		\begin{align*}
			\int_{\mathbb{R}}\left\vert E\left[ \exp (i(u\Xi +\Phi ))\right]
			\right\vert ^{2}\mathrm{d}u
			=&2\pi \int_{\mathbb{R}}(g(x)f_{\Xi }(x))^{2}\mathrm{d}x \\
			\leq &2\pi \int_{\mathbb{R}}(f_{\Xi }(x))^{2}dx=\int_{\mathbb{R}}\left\vert
			E\left[ \exp (iu\Xi )\right] \right\vert ^{2}\mathrm{d}u\text{.}
		\end{align*}%
		So%
		\begin{equation*}
			\int_{\mathbb{R}_{+}}\left\vert \phi (u)\right\vert ^{2}\mathrm{d}u\leq \int_{\mathbb{%
					R}}\left\vert E\left[ \exp (iu\Xi )\right] \right\vert ^{2}\mathrm{d}u<\infty,
		\end{equation*}
		where we have used use
		\cite[Theorem 1.1]{solomyak95}. The case $u<0$ is treated similarly.

	\end{proof}

	

	\begin{rem}
		Since Theorem~\ref{t:integrability_FT_wintnercases} holds for almost every
		$\gamma\in\left(\frac12,1\right)$, it may be viewed as extending the
		Erd\H{o}s theory of Bernoulli convolutions to the setting of
		occupation measures of drifted Takagi curves. Indeed, Solomyak's 
		almost-everywhere theorem for Bernoulli convolutions \cite{solomyak95}
		translates, through our representation of the occupation measure, into the
		existence of square-integrable local times for almost every parameter.
	\end{rem}

	Let us finally address the problem of occupation densities for the non drifted Takagi function. Consulting Corollary \ref{c:representation_simple_H}, we have to choose $\beta=1$ there, and have to face the Fourier transform of
	\begin{align}\label{e:representation_simple_H_y_0}
		\nonumber& \cT(y)-\cT(0)\\
		\nonumber =&\frac{1+\kappa}{2} S(\xi,0) - \frac{\kappa^2}{2}\frac{\gamma}{1-\gamma}\\
		&  + \sum_{n=1}^\infty X_n\, \big[(1+\kappa) S(\xi,0) (\eh)^{n+1} + \kappa^2 \gamma^n (1-\sum_{\ell=1}^\infty (\eh)^{\ell+1} X_{n+\ell})\big].
	\end{align}

	Its square, the integrability properties of which we have to assess, will not contain the terms $\frac{\kappa}{2} S(\xi,0) - \frac{\kappa^2}{2}\frac{\gamma}{1-\gamma}.$ Hence we have to deal with the Fourier transform of
	\begin{align}\label{e:representation_simple_H_essential}
		&\sum_{n=1}^\infty X_n\, \big[(1+\kappa) S(\xi,0) (\eh)^{n+1} + \kappa^2 \gamma^n (1-\sum_{\ell=1}^\infty (\eh)^{\ell+1} X_{n+\ell})\big]\\
		\nonumber =& \sum_{n=1}^\infty \gamma^n X_n\, \big[\kappa^2 + \frac{(1+\kappa)\kappa^{n}}{2} S(\xi,0)] - \frac{\kappa^2}{2} \sum_{n<m} \kappa^{-n} 2^{-m} X_n X_m\\
		\nonumber=&I_1(f_n) + I_2(\tilde{g}),
	\end{align}
	with
	$$f_n = \gamma^n [\kappa^2 + \frac{(1+\kappa)\kappa^{n}}{2} S(\xi,0)],\quad \tilde{g} = -\frac{\kappa^2}{2} g,$$
	and $g$ as chosen as in Lemma \ref{l:symmetrization:g}.

	In the sequel, We want to show for $\gamma =2^{-\frac{1}{m}},m\geq 7$  that $I_{1}(f)$ has a square integrable probability density.	To this end, define the
	functions $g_{n}$ given by%
	\begin{equation*}
		g_{k}(u):=a_{k}u\text{,}
	\end{equation*}%
	where $a_{k}:=A+C\kappa ^{k}$ with $A:=\kappa ^{2}$ and $C:=\frac{1}{2}%
	(1+\kappa )S(\xi ,0)$. Then%
	\begin{equation*}
		\left\vert \phi (u)\right\vert ^{2}=\dprod\limits_{n\geq 1}\cos
		^{2}(v\lambda ^{n}a_{n})=\dprod\limits_{q=0}^{m-1}\Psi _{q}(2^{-\frac{q}{m}%
		}v),
	\end{equation*}%
	where%
	\begin{equation*}
		\Psi _{q}(w):=\dprod\limits_{k\geq 1}\cos ^{2}(g_{m k+q}(2^{-k}w)\text{.}
	\end{equation*}
	
	In the remainder of this section, we focus on the function $\Psi_q$, whose argument we continue to denote by $v$. We now outline the strategy for proving the integrability of $\phi^2$. To this end, we decompose the infinite product defining $\Psi$ into two factors according to the size of the index $k$. The contribution of the large indices is estimated using a general argument, whereas the small indices require a more delicate analysis. The latter relies on Weyl's equidistribution theorem and suitable variants thereof; see, for example, Kuipers and Niederreiter \cite{kuipersniederreiter74}.

	To separate the two parts, for $v\ge 0$ let
	$$
	L_q(v) = \inf\{ k\ge 1: g_{m k+q}(2^{-k}v)\le \frac{\pi}{4}\},\quad q=0,\ldots, m-1.
	$$
	It is easy to see that $L_q(v)$ is well defined. Now we simply estimate the second part by using boundedness of $\cos$. So
	\be\label{e:FT-estimate_large_indices}
	\prod_{k\ge L_q(v)} \cos^2(g_{mk+q}(2^{-k}v)) \le 1.
	\ee
	
	Let us estimate $L_q(v)$ for $v\ge 0$. From the definition of $L_q(v)$, it follows that 
	
	$$
	A2^{-L_q(v)}v \le \frac{\pi}{4} < (A+C)(2^{-L_q(v)-1}v).
	$$
	So we obtain the inequalities
	\be\label{e:estimate_L(v)2}
	\log_2\big(\frac{4}{\pi}A\cdot v\big) \leq L_q(v)< \log_2\big(\frac{8(A+C)}{\pi}v\big).
	\ee
	
	We now enter into the estimation of the part of the infinite product of our Fourier transform associated with factors of small indices, namely
	\be\label{e:large_indices}
	\rho(v):= \prod_{1\le k< L_q(v)} \cos^2(g_{m k+q}(2^{-k}v)) = \exp(\sum_{1\le k< L_q(v)} \ln \cos^2 g_{m k+q}(2^{-k}v)).
	\ee
	For $v\in[2^{p-1},2^p]$ we have $(p-1) \le \log_2 (v)\le p .$ Therefore, according to \eqref{e:estimate_L(v)2}, the summation in $k$ in \eqref{e:large_indices} has to extend from $1$ to $\log_2\big(\frac{8(A+C)}{\pi}\big)+p.$ If we use the transformation $v := 2^p w$ for $w\in[\eh,1]$, we have
	$$g_{m k+q}(2^{-k}v) = g_{m k+q}(2^{p-k}w) = g_{m(p-l) +q}(2^l w),\quad l=p-k,$$
	with $\log_2\big(\frac{4A}{\pi}\big)\le l < \log_2\big(\frac{8(A+C)}{\pi}\big)+p.$ Since $\ln \cos^2$ is periodic with period $\pi$, we can define
	$$x^p_{l,q}(w) := g_{m((p-l)\vee 1)+q}(2^l w) \quad \mbox{(mod $\pi$)}),\quad l\ge 1, p\in\N,$$
	to identify
	$$ \sum_{\log_2\big(\frac{4A}{\pi}\big)\le l < \log_2\big(\frac{8(A+C)}{\pi}\big)+p} \ln \cos^2(g_{m(p-l)+q}(2^l w)) = \sum_{\log_2\big(\frac{4A}{\pi}\big)\le l < \log_2\big(\frac{8(A+C)}{\pi}\big)+p} \ln \cos^2(x^p_{l,q}(w)).$$
	To discuss the asymptotic behaviour of these sums as $p\to\infty$, we will employ equidistribution of the averaged sums of the $x_{l,q}^p(w), l\ge 1.$ Since these averages do not depend on finite fixed numbers of summands, we shall from now on study the asymptotic behaviour of the (averaged) sums
	$$\sum_{1\le l\le p} \ln \cos^2(x_{l,q}^p(w)),\quad \mbox{and associated measures}\quad \mu_p = \frac{1}{p} \sum_{l=1}^p \delta_{x^p_{l,q}(w)}.$$
	Equidistribution is derived by means of a variant of Weyl's theorem, called Koksma's general metric theorem (see Kuipers, Niederreiter \cite{kuipersniederreiter74}, p. 34).


	\begin{lem}\label{l:koksma}
		Assume that $m\geq 7$. Define for $q=0,\ldots,m-1$ and $p\geq 1$ and $%
		w\in \left[ \frac{1}{2},1\right] $%
		\begin{equation*}
			u_{l,q}^{p}(w):=g_{m((p-l)\vee 1)+q}(2^{l}w)=2^{l}wa_{m((p-l)\vee 1)+q}.
		\end{equation*}%
		Then
		\begin{equation*}
			x_{l,q}^{p}(w)=u_{l,q}^{p}(w)\text{ \ }(\textrm{mod}\,\pi )\text{.}
		\end{equation*}%
		and for all $q$ the sequence $(\mu_{p})_{p\in \mathbb{N}}$ converges to
		equidistribution on $\left[ 0,\pi \right] $ for \textup{(}Lebesgue\textup{)} almost all $w\in \left[ \frac{%
			1}{2},1\right] $.
	\end{lem}
	
	\begin{proof}
		Without loss of generality let $l_{1}>l_{2}$. Then we observe that%
		\begin{equation*}
			\frac{\mathrm d}{\mathrm dw}(u_{l_{1},q}^{p}(w)-u_{l_{2},q}^{p}(w))=2^{l_{1}}a_{m((p-l_{1})%
				\vee 1)+q}-2^{l_{2}}a_{m((p-l_{2})\vee 1)+q}
		\end{equation*}%
		is monotone in $w$. 
		Consider now the function $f:(0,\infty) \rightarrow \mathbb{R}$ given by 
		$$
		f(l):=2^{l}a_{m((p-l)\vee 1)+q}=2^{l}(A+C2^{(1-m)((p-l)\vee 1)+{\frac{1-m}{m}}q})
		$$
		The function $f$ is absolutely continuous on compact intervals in $(0,\infty)$ and its
		derivative is given
		$$
		f'(l)=
		\begin{cases}
			\ln(2)\big(A2^{l}+mC2^{(1-m)p+ml+\frac{1-m}{m}}q\big), \quad \text{ if } l<p-1\\
			\ln(2)2^{l}\big(A+mC2^{(1-m)+\frac{1-m}{m}q}\big) , \quad \text{ if } l>p-1.
		\end{cases}
		$$
		If $l<p-1$, then 
		$$
		\big|mC2^{(1-m)p+ml+\frac{1-m}{m}}q\big|\leq m|C|2^{(1-m)(p-l)}\leq m|C|2^{(1-m)}.
		$$
		On the other hand, we know that 
		$$
		\frac 12\le \kappa=2^{-1+\frac 1m}\leq \frac{1}{\sqrt{2}} \text{ as well as  } |S(\xi ,0)|\leq 1+\sqrt{2}.
		$$
		So $A=\kappa^2\geq \frac 14$ and $|C|\le \frac{(1+\sqrt{2})^2}{2\sqrt{2}}=:B$. Hence,
		$$
		f'(l)\geq \ln(2)2^{l}\big(\frac14-mB2^{-(m-1)}\big)\geq \ln(2)2^{l}0.024 \quad \text{ for all } m\geq 7.
		$$
		On the other hand, if $l>p-1$, then
		$$
		f'(l)\geq \ln(2)2^{l}\big(\frac14-B2^{-(m-1)}\big)\geq \ln(2)2^{l}\big(\frac14-mB2^{-(m-1)}\big)\geq \ln(2)2^{l}0.024 \quad \text{ for all } m\geq 7.
		$$
		Therefore 
		$$
		f'(l)\geq\ln(2)2^{l}0.024 \quad \text{ for all }  l\neq p-1, \quad \text{ and all }  m\geq 7
		$$.
		Thus we get that
		\begin{eqnarray*}
			\frac{\mathrm d}{\mathrm dw}(u_{l_{1},q}^{p}(w)-u_{l_{2},q}^{p}(w)) 
			=\int_{l_1}^{l_2}f'(l)\mathrm{d}l\geq \ln(2)0.024 \int_{l_1}^{l_2}2^{l}\mathrm{d}l  
			\geq 0.024(2^{l_{1}}-2^{l_{2}})
		\end{eqnarray*}%
		for all  $m\geq 7$.
		So the proof follows from Koksma`s general metric theorem (see \cite[p. 34]
		{kuipersniederreiter74}).
	\end{proof}

	In our study of equidistribution we can go one step ahead, and use the essential steps of the proof of Lemma \ref{l:koksma} to give an estimate for the exponential deviation from functional empirical averages and the integral of the function with respect to the limiting Lebesgue measure on $[0,\pi]$. More precisely, we will estimate for a function $h:[0,\pi]\to\R$ of bounded variation $V(h)$ the asymptotic behavior in $p$ of the integral in $w$ of
	\bean\label{e:estimate_expofaveragedsums_multiplicative}
	&&\exp(p \big[ \sum_{q=0}^{m-1} \frac{1}{p}\sum_{1\le l\le p} h(g_{m((p-l)\vee 1)+q}(2^{l} 2^{-\frac{q}{m}}w)) - m  \frac{1}{\pi}\int_0^\pi h(x) \mathrm dx\big])\\
	\nonumber&&\hspace{.5cm} \exp(p \big[ \sum_{q=0}^{m-1} \frac{1}{p} \sum_{1\le l\le p} h(x_{l,q}^p( 2^{-\frac{q}{m}}w)) - m \frac{1}{\pi}\int_0^\pi h(x) \mathrm  dx\big]).
	\eean
	This in fact will turn out to turn out to depend essentially on estimating the \textit{discrepancy} of the $x_{l,q}^p(w), 1\le l\le p.$ The discrepancy of $x_{l,q}^p(w), 1\le l\le p,$ is given by
	$$D_p := \sup_{0\le a<b\le \pi} |\frac{1}{p} \sum_{1\le l\le p} 1_{[a,b[}(\delta_{x_{l,q}^p(w)}) - (b-a)|.$$
	According to Erdös-Turan (see Kuipers, Niederreiter \cite{kuipersniederreiter74}, p. 112), we have for any $q^\ast\in\N$
	\be\label{e:erdosturan}
	D_p\le \frac{6}{q^\ast+1} + \frac{4}{\pi} \sum_{k=1}^{q^\ast} \big(\frac{1}{k} - \frac{1}{q^\ast+1}\big) \big|\frac{1}{p} \sum_{l=1}^p \exp(2\pi i x_{l,q}^p(w) k) \big|.
	\ee
	Not to overload the notation, let us state the following Lemma for the simple case of just one summand in the exponential of \eqref{e:estimate_expofaveragedsums_multiplicative}.

	\begin{lem}\label{l:BV_discrepancy}
		Let $q\in \left\{ 0,\ldots,m-1\right\} , m\geq 7\,$, $p\geq 1$ and $h:\left[ 0,\pi \right]
		\longrightarrow \mathbb{R}$ be a function of bounded variation $V(h)$. Then
		there exists a constant $c>0$ such that%
		\begin{equation*}
			\int_{\frac{1}{2}}^{1}\exp \big(p\big[ \frac{1}{p}\sum_{1\leq l\leq
				p}h(x_{l,q}^{p}(w))-\frac{1}{\pi }\int_{0}^{\pi }h(x)\mathrm dx\big] \big)\mathrm dw\leq \exp
			(cV(h)p^{\frac{1}{2}}\log (p))\text{.}
		\end{equation*}
	\end{lem}
	
	\begin{proof}
		In the following constants $c$ appearing in our inequalities are universal for the essential parameters of the expressions and may change their values from line to line. We first employ the inequality of Koksma (Kuipers, Niederreiter \cite{kuipersniederreiter74}, p. 143, and p. 91) and then the inequality of Erdös-Turan (see Kuipers, Niederreiter \cite{kuipersniederreiter74}, p. 112), to obtain for $p, q^\ast\in\N$
		\bean\label{e:exponential_erdosturan}
		&&\hspace{1cm}\int_\eh^1 \exp(p \big[ \frac{1}{p} \sum_{1\le l\le p} h(x_{l,q}^p(w)) -  \frac{1}{\pi}\int_0^\pi h(x) \mathrm dx\big]) \mathrm dw\\
		\nonumber &&\hspace{.3cm} \le \int_\eh^1\exp(V(h) p D_p)\mathrm dw\\
		\nonumber&& \le \exp(6 V(h) \frac{p}{q^\ast})\,\, \int_\eh^1 \exp\big(\frac{4}{\pi} V(h) p\big[\sum_{k=1}^{q^\ast} \frac{1}{k} \big|\frac{1}{p} \sum_{l=1}^p \exp(2\pi i  x_{l,q}^p(w) k) \big|\big]\big) \mathrm dw.
		\eean
		We continue estimating the integral in the last line of \eqref{e:exponential_erdosturan}, starting with even moments of the exponent. In fact, for $m\in\N$, using the notation of Lemma \ref{l:koksma} and $\sum_{k=1}^{q^\ast} \frac{1}{k}\le \ln(q^\ast)$, we have
		\begin{align}\label{e:multiple_koksma}
			& \big(\sum_{k=1}^{q^\ast} \frac{1}{k} \frac{1}{p}|\sum_{1\le l\le p} \exp(2\pi i x_{l,q}^p(w) k)|\big)^{2m}\\
			=&\notag \big(\sum_{k=1}^{q^\ast} \frac{1}{k} \frac{1}{p}|\sum_{1\le l\le p} \exp(2\pi i u_{l,q}^p(w) k)|\big)^{2m}\\
			\nonumber &\le \big(\ln(q^\ast)\sum_{k=1}^{q^\ast} \frac{1}{k} \frac{1}{p^2}[ p + \sum_{1 \le h < l \le p} \exp(2\pi i (u_{l,q}^p(w)-u_{h,q}^p(w)) k)]\big)^{m}\\
			\nonumber& = \ln(q^\ast)^m p^{-2m}\big( \sum_{r=0}^m \binom{m}{r} (p \ln(q^\ast))^r  (\sum_{k=1}^{q^\ast} \frac{1}{k}\sum_{1 \le h < l \le p} \exp(2\pi i (u_{l,q}^p(w)-u_{h,q}^p(w)) k))^{m-r}\big)
		\end{align}
		We continue estimating the last factor in \eqref{e:multiple_koksma}, employing again elements of the argument leading to the proof of Lemma \ref{l:koksma}. In fact, we have
		
		\bea
		&&\sum_{k=1}^{q^\ast} \frac{1}{k}\sum_{1 \le h < l \le p} \exp(2\pi i (u_{l,q}^p(w)-u_{h,q}^p(w)) k))^{m-r}\\
		\nonumber&&\hspace{.3cm} = \sum_{1\le k_1,\cdots,k_{m-r}\le q^\ast} \frac{1}{k_1\cdots k_{m-r}} \sum_{s=1}^{m-r}\sum_{1\le h_s<j_s\le p} \exp(2\pi i  \sum_{s=1}^{m-r}(u_{j_s,q}^p(w) - u_{h_s,q}^p(w)) k_s)).
		\eea
		Recall from the proof of Lemma \ref{l:koksma} the monotonicity properties of the functions $u_{j,q}^p-u_{h,q}^p$ for $h<j$ to define the monotone functions
		$$f_{h_1,j_1,\cdots, h_{m-r},j_{m-r}}(z) = \sum_{s=1}^{m-r}(u_{j_s,q}^p(w) - u_{h_s,q}^p(w)) k_s),\quad 1\le h_s<j_s\le p, 1\le s\le m-r.$$
		We next integrate \eqref{e:multiple_koksma} in $w$ over the interval $[\eh,1]$, change variables w.r.t. the monotone functions $f_{h_1,j_1,\cdots, h_{m-r},j_{m-r}}$ and use the fact that the integral of $\exp(2\pi i z)$ over bounded intervals is bounded, to estimate \eqref{e:multiple_koksma} by
		\begin{align*}
			&\sum_{1\le k_1,\cdots,k_{m-r}\le q^\ast} \frac{1}{k_1\cdots k_{m-r}}\sum_{s=1}^{m-r}  \\
			\times &\sum_{1\le h_s<j_s\le p} \int_{f_{h_1,j_1,\cdots, h_{m-r},j_{m-r}}([\eh,1])} \frac{1}{\sum_{s=1}^{m-r} ((u_{j_s,q}^p)'-u_{h_s,q}^p)')(f^{-1}_{h_1,j_1,\cdots, h_{m-r},j_{m-r}}(z))k_s} \exp(2\pi i z) \mathrm dz.
		\end{align*}
		Now we use that the integral of $\exp(2\pi i z)$ in $z$ is bounded, and that (see the proof of Lemma \ref{l:koksma}) for $1\le h_s<j_s\le p, 1\le s\le m,$ we have
		$$ \sum_{s=1}^{m-r} ((u_{j_s,q}^p)'-(u_{h_s,q}^p)')(f^{-1}_{h_1,j_1,\cdots, h_{m-r},j_{m-r}}(z))k_s\ge A^2\,\sum_{s=1}^{m-r} (2^{j_s}-2^{h_s}) = A^2\,\sum_{s=1}^{m-r} 2^{h_s}(2^{j_s-h_s}-1).$$
		
		It is easy to see that
		$$\sum_{s=1}^{m-r}\sum_{1\le h_s<j_s\le p} \frac{1}{\sum_{s=1}^{m-r} 2^{h_s}(2^{j_s-h_s}-1)k_s}$$
		is bounded in $p$. Hence, using again that $\sum_{1\le k\le q^\ast} \frac{1}{k} \le \ln(q^\ast)$, we may bound the $w$-integral of \eqref{e:multiple_koksma} by 
		\bean\label{e:multiple_koksma2}
		\nonumber	&&\int_\eh^1 (\ln(q^\ast)^m p^{-2m}\big( \sum_{r=0}^m \binom{m}{r} (p \ln(q^\ast))^r  (\sum_{k=1}^{q^\ast} \frac{1}{k}\sum_{1 \le h < l \le p} \exp(2\pi i (u_{l,q}^p(w)-u_{h,q}^p(w)) k))^{m-r}\big) \mathrm dw\\
		&&\hspace{1cm}\le c (p \ln(q^\ast)^{2m} \sum_{r=0}^m \binom{m}{r} p^r = c (p \ln(q^\ast)^{2m} (p+1)^m \le c (\frac{\ln(q^\ast)}{p^{\eh}})^{2m}.
		\eean
		We finally estimate odd moments by first employing Cauchy-Schwarz' inequality, with an analogous result. So we finally have, continuing \eqref{e:exponential_erdosturan}
		\begin{align*}
			&\exp(6 V(h) \frac{p}{q^\ast})\,\, \int_\eh^1 \exp(\frac{4}{\pi} V(h) p[\sum_{k=1}^{q^\ast} \frac{1}{k} |\frac{1}{p} \sum_{l=1}^p \exp(2\pi i  x_{l,q}^p(w) k) |]) \mathrm dw\\
			\le& \exp(6 V(h) \frac{p}{q^\ast})\,\, \sum_{m=0}^\infty \frac{1}{(2m)!} (V(h) c p)^m (\frac{\ln(q^\ast)}{p^{\eh}})^m = \exp(c V(h) [\frac{p}{q^\ast} + p^{\eh} \ln(q^\ast)]).
		\end{align*}
		Now take $q^\ast=p$ to arrive at the claimed inequality.
	\end{proof}
	
	We also have the following result
	\begin{pr}	
		\label{Equidistribution}
		Let $\gamma =2^{-\frac{1}{m}}$ for some $m\geq 7$
		and denote by $\phi (u)$ the Fourier transform of $I_{1}(f)$. Then $\phi$ is square integrable on $%
		\mathbb{R}$. Hence, $I_{1}(f)$ has a square integrable probability
		density.
	\end{pr}
	\begin{proof}
		We have to estimate the integral of $v\mapsto  \phi^2(v)$ on $\R_+$, remarking that the evenness of $\cos$ trivially extends this to all of $\R.$ We obtain, with a constant $C$ standing for the contribution of large indices in the infinite product describing the Fourier transform, and a constant $p_0$ marking the domain for which small index contributions become relevant
		\bean\label{e:final_estimate_integral_FT}
		&&\hspace{.2cm}  \int_0^\infty \phi^2(v) \mathrm dv\\
		\nonumber&&\leq  C + \sum_{p = p_0}^\infty \int_{2^{p-1}}^{2^p}  \phi^2(v) dv\\
		\nonumber&& = C + \sum_{p=p_0}^\infty \int_{\eh}^1 2^{p}  [\exp(p \sum_{q=0}^{m-1} [\frac{1}{p} \sum_{1\le l\le p} \ln \cos^2(x_{l,p}^p( 2^{-\frac{q}{m}}w))])  \mathrm dw\\
		\nonumber&& = C + \sum_{p=p_0}^\infty \int_{\eh}^1\big\{  \exp\big(p \sum_{q=0}^{m-1} [\frac{1}{p} \sum_{1\le l\le p} \ln \cos^2(x_{l,p}^p( 2^{-\frac{q}{m}}w))] - m p \frac{1}{\pi} \int_0^{\pi} \ln \cos^2(x) \mathrm  dx\big)\\
		\nonumber&&\hspace{2cm} \times\exp\big[p\ln(2)+ m p \frac{1}{\pi} \int_0^{\pi} \ln \cos^2(x)  \mathrm dx\big]\big\}  \mathrm dw.
		\eean
		The Lebesgue integral of the $\pi$-periodic $\ln \cos^2$ is well known. We have
		\be\label{e:equidistribution_convergence_lncos}
		\frac{1}{\pi}\int_0^\pi \ln \cos^2(x)  \mathrm dx = \frac{4}{\pi} \int_0^{\frac{\pi}{2}} \ln \cos(x)  \mathrm dx = - 2 \ln(2).
		\ee
		Hence the last line of \eqref{e:final_estimate_integral_FT} simply becomes
		\be\label{e:FT_convergence_generating}
		\exp(p\ln(2)[1-2m]).
		\ee
		We have to be a little careful with estimating the contributions of the first line on the rhs of \eqref{e:final_estimate_integral_FT} by means of Lemma \ref{l:BV_discrepancy}, since $\ln \cos^2$ is not of bounded variation on $[0,\pi].$ Given $\epsilon>0$, choose a function of bounded variation $h:[0,\pi]\to \R$ such that $\frac{1}{\pi}\int |h(x)- \ln \cos^2(x)| \mathrm dx < \epsilon.$ Using $h$ instead of $\ln\cos^2$ in \eqref{e:final_estimate_integral_FT} forces us to modify \eqref{e:FT_convergence_generating} to
		\be\label{e:FT_convergence_generating2}
		\exp(p\ln(2)[1-2(m-\epsilon)]).
		\ee
		We apply Lemma \ref{l:BV_discrepancy} with $h$ to obtain an estimate for the contributions in the first line on the rhs of \eqref{e:final_estimate_integral_FT} which is given by
		\be\label{e:multiple_koksma3}
		\exp(c V(h) p^\eh \ln(p)).
		\ee
		But the product of \eqref{e:FT_convergence_generating2} and \eqref{e:multiple_koksma3} remains summable, since $p^\eh \ln(p)$ is dominated by $p$ as $p\to\infty.$ This proves the result. 
	\end{proof}
	
	We can finally state our main result.
	\begin{thm}
		Let $\gamma =2^{-\frac{1}{m}}$ for some $m\geq 7$. Then $\phi \in L^{2}(\mathbb{R})$. Therefore $%
		(y\longmapsto \mathcal{T}(y))$ has a square integrable occupation density.
	\end{thm}
	
	\begin{proof}
		We can proceed in the same way as in the proof of Theorem \ref{t:integrability_FT_wintnercases}
		\begin{equation}
			\phi (u)\leq \Big\vert E\Big[ \exp (i(uI_{1}(f)+\left( u\right)
			^{1/2}I_{2}(\widetilde{g})))\Big] \Big\vert   \notag
		\end{equation}%
		for all $u>0$, where%
		\begin{equation*}
			f(n)=\gamma ^{n}\big[ \kappa ^{2}+\frac{(1+\kappa )\kappa ^{n}}{2}S(\xi ,0)%
			\big] \text{ and }\widetilde{g}=-\frac{\kappa ^{2}}{2}g
		\end{equation*}%
		with $g$ as chosen before. By iterating the latter inequality and by
		applying Plancherel`s Theorem and Proposition \ref{Equidistribution} we
		similarly get%
		\begin{equation*}
			\int_{\mathbb{R}_{+}}\big\vert \phi (u)\big\vert ^{2}\mathrm du\leq \int_{\mathbb{%
					R}}\big\vert E\big[ \exp (iuI_{1}(f))\big] \big\vert ^{2}\mathrm du<\infty 
			\text{.}
		\end{equation*}%
		The case $u<0$ is treated similarly.
	\end{proof}


\end{document}